\documentclass[twoside,b5paper,11pt,english,openany]{book}
\usepackage[top=3.35cm, left=2.3cm, bottom=3.35cm, right=2.3cm,headsep=.63cm]{geometry}

\usepackage{stix}    
\pdfmapfile{=stix.map}

\usepackage{tgtermes}
\pdfmapfile{=qtm.map}
\usepackage{tgheros}
\pdfmapfile{=qhv.map}
\usepackage{tgcursor}
\pdfmapfile{=qcr.map}

\makeatletter
\newcommand\incl[3]{
  \@ifclassloaded{#1}{#3}{#2}
}
\newcommand\icl[3]{
  \@ifclassloaded{#1}{#2}{#3}
}
\newcommand\inpl[3]{
  \@ifpackageloaded{#1}{#3}{#2}
}
\newcommand\ipl[3]{
  \@ifpackageloaded{#1}{#2}{#3}
}
\makeatother

\makeatletter
\newcommand*\icBook{\includegraphics[height=\f@size pt]{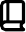}}
\newcommand*\icBookB{\includegraphics[height=\f@size pt]{material-symbols--book-2-outline-rounded~beamer.pdf}}
\newcommand*\icArticle{\includegraphics[height=\f@size pt]{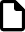}}
\newcommand*\icArticleB{\includegraphics[height=\f@size pt]{material-symbols--docs-outline-rounded~beamer.pdf}}
\newcommand*\icInbook{\includegraphics[height=\f@size pt]{material-symbols--book-ribbon-outline.pdf}}
\newcommand*\icInbookB{\includegraphics[height=\f@size pt]{material-symbols--book-ribbon-outline~beamer.pdf}}
\makeatother

\usepackage[utf8]{inputenc}
\usepackage[T1]{fontenc}
\usepackage{textcomp}
\usepackage{babel}

\usepackage[normalem]{ulem}
\usepackage{setspace}
\usepackage{color}
\incl{beamer}{\usepackage[dvipsnames]{xcolor}}{}
\usepackage{graphicx}
\incl{beamer}{\usepackage[bookmarksdepth=subsection]{hyperref}}{}
\usepackage[all]{nowidow}
\usepackage[format=plain,justification=centering]{caption}

\usepackage{url}
\usepackage{wrapfig}
\usepackage{subcaption}
\usepackage{dashbox}
\usepackage{pifont}
\usepackage{xhfill}
\incl{beamer}{\usepackage{titling}}{}

\usepackage{xargs}
\usepackage{suffix}
\usepackage{calc}
\usepackage{trimspaces}
\usepackage{lipsum}
\usepackage{etoolbox}
\usepackage{float}
\usepackage{afterpage}
\usepackage{csquotes}

\usepackage{datetime}
\newdateformat{monthyeardate}{\monthname[\THEMONTH] \THEYEAR}
\newdateformat{monthyeardaydate}{\monthname[\THEMONTH] \THEDAY, \THEYEAR}
\newdateformat{daymonthyeardate}{\THEDAY\ \monthname[\THEMONTH] \THEYEAR}
\newdateformat{defaultdateformat}{\iflanguage{french}{\THEDAY\ \monthname[\THEMONTH] \THEYEAR}{\monthname[\THEMONTH] \THEDAY, \THEYEAR}}

\incl{beamer}{
  \usepackage[explicit,compact]{titlesec}
  \titleformat{\chapter}[block]
  {\centering}{\large \textsc{\chaptername}\ $\thechapter$\\}{0ex}{\LARGE #1}
  \titlespacing*{\chapter}{0pt}{-1em}{2em}
  \titleformat{\section}
  {\bfseries\Large}{\thesection}{1em}{#1}
}{}

\icl{book}{
  \usepackage[backend=biber,style=alphabetic,sorting=nyt,url=false,isbn=false,refsection=chapter]{biblatex}
}{
  \usepackage[backend=biber,style=alphabetic,sorting=nyt,url=false,isbn=false]{biblatex}
}

\DeclareLabelalphaTemplate{
  \labelelement{
    \field[final]{shorthand}
  }
  \labelelement{
\field[strwidth=3,strside=left,ifnames=1]{labelname}
    \field[strwidth=1,strside=left]{labelname}
  }
  \labelelement{
    \field{year}
  }
}

\renewbibmacro{in:}{}

\incl{beamer}{\renewbibmacro{begentry}{\ifentrytype{book}{{\scriptsize\icBook}\ \hspace{.5pt}}{}\ifentrytype{article}{{\scriptsize\icArticle}\ \hspace{.5pt}}{}\ifentrytype{inbook}{{\scriptsize\icInbook}\ \hspace{.5pt}}{}}}{\renewbibmacro{begentry}{\ifentrytype{book}{{\scriptsize\icBookB}\ \hspace{.5pt}}{}\ifentrytype{article}{{\scriptsize\icArticleB}\ \hspace{.5pt}}{}\ifentrytype{inbook}{{\scriptsize\icInbookB}\ \hspace{.5pt}}{}}}

\newcommandx\prettytitle[4][1=,2=\thetitle,3=\theauthor,4=\defaultdateformat\today]{
    \begin{center}
      \phantom{}\vspace{-2em}

      {\large\phantom{}\xrfill[0.175\baselineskip]{0.4pt}\ifblank{#1}{}{\ #1\ \xrfill[0.175\baselineskip]{0.4pt}}\\[1em]}

      {\LARGE#2}
      \vspace{1em}

      #3
      \vspace{1em}

      #4\\[1ex]
      \phantom{}\xrfill[0.175\baselineskip]{0.4pt}
  \end{center}
}

\newcommandx\prettytitlenoauth[3][1=,2=\thetitle,3=\defaultdateformat\today]{
    \begin{center}
      \phantom{}\vspace{-2em}

      {\large\phantom{}\xrfill[0.175\baselineskip]{0.4pt}\ifblank{#1}{}{\ #1\ \xrfill[0.175\baselineskip]{0.4pt}}\\[1em]}

      {\LARGE#2}
      \vspace{1em}

      #3\\[1ex]
      \phantom{}\xrfill[0.175\baselineskip]{0.4pt}
  \end{center}
}

\newcommand\printrefs{
  \begin{center}
    \vspace{1em}
    \scshape\Large \iflanguage{french}{Références}{References}
    \vspace{1ex}
  \end{center}
  \printbibliography[heading=none]
  \addcontentsline{toc}{section}{\iflanguage{french}{Références}{References}}
}

\WithSuffix\newcommand\printrefs*{
  \begin{center}
    \vspace{1em}
    \scshape\Large \iflanguage{french}{Références}{References}
    \vspace{1ex}
  \end{center}
  \printbibliography[heading=none]
}

\newcommand\apx{
  \appendix
  \renewcommand\chaptername{\iflanguage{french}{Annexe}{Appendix}}
  \addtocontents{toc}{\protect\setcounter{tocdepth}{0}}
}

\newcommand{\ZeroRoman}[1]{\ifcase\value{#1}\relax 0\else \Roman{#1}\fi}\renewcommand\thesection{\ZeroRoman{section}}

\author{Nathan \textsc{Benichou}}

\usepackage{expl3}
\ExplSyntaxOn
\cs_new_eq:NN \Repeat \prg_replicate:nn
\ExplSyntaxOff

\makeatletter
\renewcommandx\hrulefill[1][1=.4pt]{\leavevmode\leaders\hrule height #1\hfill\kern\z@}
\makeatother \incl{beamer}{
  \usepackage{enumitem}

\newlist{props}{enumerate}{1}
  \setlist[props]{left=0pt,align=parleft,labelsep=0em,label=($\roman*$),font=\normalfont}
  \newlist{props2}{enumerate}{1}
  \setlist[props2]{left=0pt,align=parleft,labelsep=-1ex,label=($\roman*$),font=\normalfont}
\setlist[itemize]{left=0pt,label=--} \setlist[enumerate]{left=0pt,font=\normalfont}
  \newlist{enumtimes}{enumerate}{1}
  \setlist[enumtimes]{left=0pt,align=parleft,font=\normalfont,label=($\alph*$),itemsep=0ex}\newlist{enumtimesar}{enumerate}{1}
  \setlist[enumtimesar]{left=0pt,align=parleft,font=\normalfont,label=($\arabic*$),itemsep=0ex}\newlist{enumtimesrom}{enumerate}{1}
  \setlist[enumtimesrom]{left=0pt,align=parleft,font=\normalfont,label=({\itshape\roman*}\hspace{.5pt}),itemsep=0ex}\newlist{enumtimesromp}{enumerate}{1}
  \setlist[enumtimesromp]{left=0pt,align=parleft,font=\normalfont,label=({$\textit{\roman*}'$}\hspace{.5pt}),itemsep=0ex}}{} \incl{beamer}{\usepackage{minitoc}
\usepackage{tocloft}
\usepackage{scrextend}   

\setcounter{tocdepth}{3}

\makeatletter
\newcommand*{\toccontents}{\csname @starttoc\endcsname{toc}}
\makeatother

\newcommand\printtoc{
  \begin{center}
    \vspace{1em}
    \scshape\Large \iflanguage{french}{Table des matières}{Contents}
    \vspace{1ex}
  \end{center}
  \toccontents
}

\newcommand\addssdesc[1]{
  \addtocontents{toc}{
    \begin{addmargin}[\cftsubsecindent-1.325ex]{\cftsubsecnumwidth-1.325ex}\footnotesize\raggedright #1
    \end{addmargin}
  }
} }{}
\definecolor{alizarin}{rgb}{0.905882353,0.298039216,0.235294118}
\definecolor{nephritis}{rgb}{0.152941176,0.682352941,0.376470588}
\definecolor{belizehole}{rgb}{0.160784314,0.501960784,0.725490196}
\definecolor{peterriver}{rgb}{0.203921569,0.596078431,0.858823529}
\definecolor{silver}{rgb}{0.74117647058823529412,0.76470588235294117647,0.78039215686274509804}
\definecolor{silverdk}{rgb}{0.66705882352941176471,0.68823529411764705882,0.70235294117647058824}

\definecolor{officeblueint}{rgb}{0.51372549,0.745098039,0.925490196}
\definecolor{officeblueext}{rgb}{0,0.388235294,0.694117647}
\definecolor{officegreenint}{rgb}{0.631372549,0.866666667,0.666666667}
\definecolor{officegreenext}{rgb}{0.188235294,0.564705882,0.282352941}
\definecolor{officeredint}{rgb}{1,0.568627451,0.596078431}
\definecolor{officeredext}{rgb}{0.831372549,0.137254902,0.078431373}
\definecolor{officeyellowint}{rgb}{0.97254902,0.858823529,0.560784314}
\definecolor{officeyellowext}{rgb}{0.870588235,0.423529412,0}
\definecolor{officepurpleint}{rgb}{0.831372549,0.57254902,0.847058824}
\definecolor{officepurpleext}{rgb}{0.62745098,0.294117647,0.658823529}
\definecolor{officegrayint}{rgb}{0.784313725,0.776470588,0.768627451}
\definecolor{officegrayext}{rgb}{0.474509804,0.466666667,0.454901961}
\definecolor{officeblackint}{rgb}{0.784313725,0.776470588,0.768627451}
\definecolor{officeblackext}{rgb}{0.22745098,0.22745098,0.219607843} \usepackage{tabularray}

\UseTblrLibrary{diagbox}

\newcommand\matr[1]{\left(\begin{tblr}{cells=c,rowsep=1pt,colsep=4pt}
  #1
\end{tblr}\right)}
\WithSuffix\newcommand\matr*[1]{\left(\begin{tblr}{cells=c,rowsep=0pt,colsep=3pt}
  #1
\end{tblr}\right)}

\newcommand\detr[1]{\,\begin{tblr}{cells=c,rowsep=.5pt,colsep=4pt,vline{1,Z}={.6pt}}
  #1
\end{tblr}}

\WithSuffix\newcommand\detr*[1]{\,\begin{tblr}{cells=c,rowsep=.5pt,colsep=1.5pt,vline{1,Z}={.6pt}}
  #1
\end{tblr}}

\newcolumntype{O}{>{{}}c<{{}}} \newcolumntype{C}{>{\displaystyle}c} \newcolumntype{L}{>{\displaystyle}l} \newcolumntype{R}{>{\displaystyle}r} %
 \usepackage{pgf}
\usepackage{tikz}
\usepackage{pgfplots}
\usepackage{tikz-cd}

\usetikzlibrary{arrows, arrows.meta, bending, calc, backgrounds, patterns}
\usetikzlibrary{decorations.markings}
\usetikzlibrary{bbox}
\usetikzlibrary{decorations.pathreplacing}

\usetikzlibrary {patterns.meta}
\tikzdeclarepattern{
  name=mylines,
  parameters={
      \pgfkeysvalueof{/pgf/pattern keys/size},
      \pgfkeysvalueof{/pgf/pattern keys/angle},
      \pgfkeysvalueof{/pgf/pattern keys/line width},
    },
  bounding box={
      (0,-0.5*\pgfkeysvalueof{/pgf/pattern keys/line width}) and
      (\pgfkeysvalueof{/pgf/pattern keys/size},
      0.5*\pgfkeysvalueof{/pgf/pattern keys/line width})},
  tile size={(\pgfkeysvalueof{/pgf/pattern keys/size},
      \pgfkeysvalueof{/pgf/pattern keys/size})},
  tile transformation={rotate=\pgfkeysvalueof{/pgf/pattern keys/angle}},
  defaults={
      size/.initial=5pt,
      angle/.initial=45,
      line width/.initial=.4pt,
    },
  code={
      \draw [line width=\pgfkeysvalueof{/pgf/pattern keys/line width}]
      (0,0) -- (\pgfkeysvalueof{/pgf/pattern keys/size},0);
    },
}

\newcommand*\circled[1]{\tikz[baseline=(char.base)]{
    \node[shape=circle,draw,inner sep=1pt,line width=.6pt] (char) {#1};}}

 \usepackage{amsmath}
\usepackage{amsthm}
\usepackage{mathrsfs}     \usepackage{mathtools}
\usepackage{extarrows}  
  \usepackage{marvosym}     \usepackage{bm}

\allowdisplaybreaks

\newcommand\bskip{\vspace{.5\bigskipamount}}

\makeatletter
\ipl{pgf}{
  \def\slashedarrowfill@#1#2#3#4#5{$\m@th\thickmuskip0mu\medmuskip\thickmuskip\thinmuskip\thickmuskip
     \relax#5#1\mkern-7mu\cleaders\hbox{$#5\mkern-2mu#2\mkern-2mu$}\hfill
     \mathclap{#3}\mathclap{#2}\cleaders\hbox{$#5\mkern-2mu#2\mkern-2mu$}\hfill
     \mkern-7mu#4$}

  \def\rightslashedarrowfilla@{\slashedarrowfill@\relbar\relbar{\raisebox{1.2pt}{$\scriptscriptstyle\diagup$}}\rightarrow}
  \newcommand\xslashedrightarrowa[2][]{\ext@arrow 0055{\rightslashedarrowfilla@}{#1}{#2}}

  \def\rightslashedarrowfillb@{\slashedarrowfill@\relbar\relbar/\rightarrow}
  \newcommand\xslashedrightarrowb[2][]{\ext@arrow 0055{\rightslashedarrowfillb@}{#1}{#2}}

  \def\rightslashedarrowfille@{\slashedarrowfill@\relbar\relbar{\raisebox{0em}{\footnotesize/}}\rightarrow}
  \newcommand\xslashedrightarrowe[2][]{\ext@arrow 0055{\rightslashedarrowfille@}{#1}{#2}}

  \def\rightslashedarrowfillc@{\slashedarrowfill@\relbar\relbar{\raisebox{.12em}{\tiny/}}\rightarrow}
  \newcommand\xslashedrightarrowc[2][]{\ext@arrow 0055{\rightslashedarrowfillc@}{#1}{#2}}

  \pgfdeclareshape{slash underlined}
  {
    \inheritsavedanchors[from=rectangle] \inheritanchorborder[from=rectangle]
    \inheritanchor[from=rectangle]{north}
    \inheritanchor[from=rectangle]{north west}
    \inheritanchor[from=rectangle]{north east}
    \inheritanchor[from=rectangle]{center}
    \inheritanchor[from=rectangle]{west}
    \inheritanchor[from=rectangle]{east}
    \inheritanchor[from=rectangle]{mid}
    \inheritanchor[from=rectangle]{mid west}
    \inheritanchor[from=rectangle]{mid east}
    \inheritanchor[from=rectangle]{base}
    \inheritanchor[from=rectangle]{base west}
    \inheritanchor[from=rectangle]{base east}
    \inheritanchor[from=rectangle]{south}
    \inheritanchor[from=rectangle]{south west}
    \inheritanchor[from=rectangle]{south east}
    \inheritanchorborder[from=rectangle]
    \foregroundpath{
\southwest \pgf@xa=\pgf@x \pgf@ya=\pgf@y
      \northeast \pgf@xb=\pgf@x \pgf@yb=\pgf@y
      \pgf@xc=\pgf@xa
      \advance\pgf@xc by .5\pgf@xb
      \pgf@yc=\pgf@ya
      \advance\pgf@xc by -1.3pt
      \advance\pgf@yc by -1.8pt
      \pgfpathmoveto{\pgfqpoint{\pgf@xc}{\pgf@yc}}
      \advance\pgf@xc by  2.6pt
      \advance\pgf@yc by  3.6pt
      \pgfpathlineto{\pgfqpoint{\pgf@xc}{\pgf@yc}}
      \pgfpathmoveto{\pgfqpoint{\pgf@xa}{\pgf@ya}}
      \pgfpathlineto{\pgfqpoint{\pgf@xb}{\pgf@ya}}
   }
  }
  \tikzset{nomorepostaction/.code=\let\tikz@postactions\pgfutil@empty}

  }{}
\makeatother

\newcommand\xxrightarrow[2][]{
    \xrightarrow[{\raisebox{1.25ex-\heightof{$\scriptstyle#1$}}[0pt]{$\scriptstyle#1$}}]{#2}}
\newcommand\xxlongequal[2][]{
    \xlongequal[{\raisebox{1.25ex-\heightof{$\scriptstyle#1$}}[0pt]{$\scriptstyle#1$}}]{#2}}

\makeatletter
\newlength{\osetlength}
\newcommand{\oset}[3][0ex]{\setlength{\osetlength}{#1}
  \mathrel{\mathop{#3}\limits^{
    \vbox to\osetlength{\kern-2\ex@
    \hbox{$\scriptstyle#2$}\vss}}}}
\makeatother

\newcommandx\cv[4][1=,2=,3=0ex,4=0ex]{\oset[-.25ex+#3]{#2}{\oset[-2ex+#4]{#1}{\longrightarrow}}}
\newcommandx\cvx[2][1=,2=]{\xxrightarrow[#1]{#2}}
\ipl{pgf}{\newcommandx\notcvx[2][1=,2=]{\xslashedrightarrowc[#1]{#2}}}{}
\WithSuffix\newcommandx\cvx*[2][1=,2=]{\xrightarrow[#1]{#2}}
\newcommandx\cvcs[1][1=]{\cv[#1][\textnormal{c.s.}][-0.5ex]}
\newcommandx\cvps[1][1=]{\cv[#1][\textnormal{p.s.}][-0.25ex]}
\newcommandx\cvas[1][1=]{\cv[#1][\textnormal{a.s.}][-0.25ex]}
\newcommandx\cvpp[1][1=]{\cv[#1][\textnormal{p.p.}][-0.25ex]}
\newcommandx\cvcsx[1][1=]{\cvx[#1][\textnormal{c.s.}]}
\newcommandx\cvpsx[1][1=]{\cvx[#1][\textnormal{p.s.}]}
\newcommandx\cvasx[1][1=]{\cvx[#1][\textnormal{a.s.}]}
\newcommandx\cvppx[1][1=]{\cvx[#1][\textnormal{p.p.}]}

\newcommandx\egv[4][1=,2=,3=0ex,4=0ex]{\oset[-.25ex+#3]{#2}{\oset[-2ex+#4]{#1}{=}}}
\newcommandx\egvx[2][1=,2=]{\xxlongequal[#1]{#2}}
\WithSuffix\newcommandx\egvx*[2][1=,2=]{\xlongequal[#1]{#2}}
\newcommandx\egvcsx[1][1=]{\egvx[#1][\textnormal{c.s.}]}
\newcommandx\egvpsx[1][1=]{\egvx[#1][\textnormal{p.s.}]}
\newcommandx\egvppx[1][1=]{\egvx[#1][\textnormal{p.p.}]}
\newcommandx\egvcs[1][1=]{\egv[#1][\textnormal{c.s.}][-0.75ex]}
\newcommandx\egvps[1][1=]{\egv[#1][\textnormal{p.s.}][-0.25ex]}
\newcommandx\egvpp[1][1=]{\egv[#1][\textnormal{p.p.}][-0.25ex]}

\newcommandx\eqv[4][1=,2=,3=0ex,4=0ex]{\oset[1.75ex+#3]{#2}{\oset[-2ex+#4]{#1}{\sim}}}
\newcommandx\eqvcs[1][1=]{\eqv[#1][\textnormal{c.s.}][-0.75ex]}
\newcommandx\eqvps[1][1=]{\eqv[#1][\textnormal{p.s.}][-0.25ex]}
\newcommandx\eqvas[1][1=]{\eqv[#1][\textnormal{a.s.}][-0.25ex]}
\newcommandx\eqvpp[1][1=]{\eqv[#1][\textnormal{p.p.}][-0.25ex]}

\makeatletter
\def\smallunderbrace#1{\mathop{\vtop{\m@th\ialign{##\crcr
   $\hfil\displaystyle{#1}\hfil$\crcr
   \noalign{\kern3\p@\nointerlineskip}\tiny\upbracefill\crcr\noalign{\kern3\p@}}}}\limits}
\makeatother

\newcommand\qtx[1]{\quad\text{#1}\quad}

\renewcommand\lq{\leqslant}
\newcommand\gq{\geqslant}

\newcommand\cp{^{\textnormal{\textsf{c}}}}
\newcommand\ssm{\!\smallsetminus\!}

 \DeclareMathOperator\B{\mathcal{B}}

\newcommand\sac{\kern1pt|\kern1pt}

\newcommand\N{\textnormal{\textbf{N}}} 

\newcommand\sk{\textnormal{\textbf{D}}} 
\newcommand\tstc{\mathcal{E}}
  \newcommand\D{\textnormal{d}}

\newcommand\T{\mathcal{T}}
 
 \newcommand\G{\mathcal{G}}

\newcommand\ct{\mathcal{C}}

\newcommandx\cl[1][1=\alpha]{\mathcal{H}^{#1}}
\newcommandx\clb[1][1=\alpha]{\mathcal{H}^{#1}_{\textnormal{\textsf{b}}}}
\newcommandx\clloc[1][1=\alpha]{\mathcal{H}^{#1}_{\textnormal{loc}}}

\DeclareMathOperator{\id}{id}

\DeclareMathOperator{\dist}{dist}

\newcommand\R{\textnormal{\textbf{R}}}  
 \newcommand\Q{\textnormal{\textbf{Q}}} 
\newcommand\F{\mathcal{F}}

\newcommand\lpl{L_{\textnormal{loc}}}
\WithSuffix\newcommand\lpl*{L_{\textnormal{(loc)}}}
\newcommand\pt{\mathcal{P}}

\DeclareMathOperator\lip{Lip}
\newcommand\hol{H}

\DeclareMathOperator\sym{\mathfrak{S}}
\incl{beamer}{}{}
\newcommandx\W[2][1=1,2=1]{W^{#2,#1}}
\newcommandx\wl[2][1=1,2=1]{W_{\textnormal{loc}}^{#2,#1}}
\newcommandx\hml[1][1=1]{H_{\textnormal{loc}}^{#1}}
\WithSuffix\newcommandx\wl*[2][1=1,2=1]{W_{\textnormal{(loc)}}^{#2,#1}}
\newcommand{\interior}[1]{{\kern0pt#1}^{\circ}}

\newcommand\ind{\textnormal{\textbf{1}}}

\DeclareMathOperator*{\Li}{Li}
\DeclareMathOperator*{\Ls}{Ls}

\DeclareMathOperator{\osc}{osc}

\newcommand\sm[1]{\left\llbracket #1\right\rrbracket}
\WithSuffix\newcommand\sm*[1]{\llbracket #1\rrbracket}
\newcommand\sms[1]{\left(\!\left(#1\right)\!\right)}
\WithSuffix\newcommand\sms*[1]{(\!(#1)\!)}
\WithSuffix\newcommand\sqrt*[1]{\sqrt{\phantom{|}\!\!\smash{#1}}}
\newcommand\abs[1]{|#1|}

\WithSuffix\newcommand\abs*[1]{\left|#1\right|}
\newcommand\norm[1]{\|#1\|}

\newcommand\normd{\norm{\kern1pt\cdot\kern1pt}}
\WithSuffix\newcommand\norm*[1]{\left\|#1\right\|}
\newcommand\fun[5]{#1:\begin{aligned}
    #2 \longrightarrow & \ #3 \\
    #4    \longmapsto     & \ #5
\end{aligned}}
\WithSuffix\newcommand\fun*[4]{\begin{aligned}
    #1 \longrightarrow & \ #2 \\
    #3    \longmapsto     & \ #4
\end{aligned}}
\newcommandx\matelem[3][1=,2=]{\phantom{}^{#1}[#3]_{#2}}

\DeclareMathOperator*{\adh}{adh}
\newcommand\sys[1]{\left\{\begin{aligned}#1\end{aligned}\right.}
\WithSuffix\newcommand\sys*[1]{\begin{aligned}#1\end{aligned}}

\newcommand\bigzero{\makebox(0,0){\textnormal{\fontfamily{lmr}\selectfont\huge0}}}
\WithSuffix\newcommand\bigzero*{\textnormal{\fontfamily{lmr}\selectfont\huge0}}

\newcommandx\psb[3][1=\middle|]{\left\langle \left.{#2} \kern1.25pt#1\kern1.25pt {#3}\right. \right\rangle}
\newcommand\pss[1]{\left\langle \left.{#1}\right. \right\rangle}
\WithSuffix\newcommand\pss*[1]{\langle \left.{#1}\right. \rangle}

\newcommand\ps[2]{\left\langle \left.{#1} \kern1.25pt\middle|\kern1.25pt {#2}\right. \right\rangle}

\WithSuffix\newcommand\ps*[2]{\langle \left.{#1} \kern1.25pt|\kern1.25pt {#2}\right. \rangle}
\newcommand\phl[2]{\stretchleftright{\langle\kern-1.59pt\rule[-.2pt]{.45pt}{5.5pt}}{\kern-1.1pt\left.{#1} \kern1.25pt\middle|\kern1.25pt {#2}\right.\kern-3pt}{\rangle}\,}
\newcommand\phr[2]{\stretchleftright{\langle}{\kern-3pt\left.{#1} \kern1.25pt\middle|\kern1.25pt {#2}\right.\kern-2pt}{\rangle\kern-2.5pt\rule[-.2pt]{.45pt}{5.5pt}}\,}

\newcommand\lightqed{\textnormal{/\kern-.75pt/}}

\newtheoremstyle{noswp}{}{}{\itshape}{}{\bfseries\boldmath}{.\,}{ }{\thmname{#1}\thmnumber{ #2}\thmnote{ (#3)}}

\newtheoremstyle{noswpdefn}{}{}{\normalfont}{}{\bfseries\boldmath}{.\,}{ }{\thmname{#1}\thmnumber{ #2}\thmnote{ (#3)}}

\newtheoremstyle{swp}{}{}{\itshape}{}{\bfseries\boldmath}{.\,}{ }{\thmnumber{#2}\thmname{ #1}\thmnote{ (#3)}}
  
\newtheoremstyle{swpt}{}{}{\itshape}{}{\bfseries\boldmath}{.\,}{ }{\thmnumber{#2.}\thmname{ #1}\thmnote{ (#3)}}

\newtheoremstyle{swpdefn}{}{}{\normalfont}{}{\bfseries\boldmath}{.\,}{ }{\thmnumber{#2}\thmname{ #1}\thmnote{ (#3)}}
  
\newtheoremstyle{swptdefn}{}{}{\normalfont}{}{\bfseries\boldmath}{.\,}{ }{\thmnumber{#2.}\thmname{ #1}\thmnote{ (#3)}}

\makeatletter
\def\th@plain{\thm@notefont{}\itshape }
\def\th@definition{\thm@notefont{}\normalfont }
\makeatother

\newtheorem{innercustomthm}{\textsc{\iflanguage{french}{Théorème}{Theorem}}}
\newenvironment{cthm}[1]
  {\renewcommand\theinnercustomthm{#1}\innercustomthm}
  {\endinnercustomthm}

\theoremstyle{swp}

\newtheorem{thm}{\textsc{\iflanguage{french}{Théorème}{Theorem}}}

\theoremstyle{noswp}

\newtheorem{nthm}[thm]{\textsc{\iflanguage{french}{Théorème}{Theorem}}}

\newtheorem{npropn}[thm]{\textsc{\iflanguage{french}{Proposition}{Proposition}}}
\newtheorem{nlem}[thm]{\textsc{\iflanguage{french}{Lemme}{Lemma}}}

\theoremstyle{plain}

\newtheorem*{thm*}{\textsc{\iflanguage{french}{Théorème}{Theorem}}}
\newtheorem*{prop*}{\textsc{\iflanguage{french}{Propriété}{Property}}}
\newtheorem*{propn*}{\textsc{\iflanguage{french}{Proposition}{Proposition}}}
\newtheorem*{lem*}{\textsc{\iflanguage{french}{Lemme}{Lemma}}}
\newtheorem*{cor*}{\textsc{\iflanguage{french}{Corollaire}{Corollary}}}
\newtheorem*{thdefn*}{\textsc{\iflanguage{french}{Théorème-Définition}{Theorem-Definition}}}
\newtheorem*{pronotatn*}{\textsc{\iflanguage{french}{Propriété-Notation}{Property-Notation}}}
\newtheorem*{propdefn*}{\textsc{\iflanguage{french}{Propriété-Définition}{Property-Definition}}}
\newtheorem*{cordefn*}{\textsc{\iflanguage{french}{Corollaire-Définition}{Corollary-Definition}}}
\newtheorem*{propndefn*}{\textsc{\iflanguage{french}{Proposition-Définition}{Proposition-Definition}}}

\theoremstyle{swpt}

\theoremstyle{remark}

\theoremstyle{swpdefn}

\theoremstyle{definition}

\newtheorem*{cadr*}{\textsc{\iflanguage{french}{Cadre}{Framework}}}
\newtheorem*{conv*}{\textsc{\iflanguage{french}{Convention}{Convention}}}
\newtheorem*{term*}{\textsc{\iflanguage{french}{Terminologie}{Terminology}}}
\newtheorem*{defn*}{\textsc{\iflanguage{french}{Définition}{Definition}}}
\newtheorem*{ex*}{\textsc{\iflanguage{french}{Exemple}{Example}}}
\newtheorem*{rk*}{\textsc{\iflanguage{french}{Remarque}{Remark}}}
\newtheorem*{rap*}{\textsc{\iflanguage{french}{Rappel}{Reminder}}}
\newtheorem*{notat*}{\textsc{\iflanguage{french}{Notation}{Notation}}}
\newtheorem*{hyp*}{\textsc{\iflanguage{french}{Hypothèse}{Assumption}}}
\newtheorem*{idprf*}{\textsc{\iflanguage{french}{Idée de la preuve}{Idea of the proof}}}

\theoremstyle{swptdefn}

\theoremstyle{noswpdefn} 

\newtheorem{nncadr}[thm]{\textsc{\iflanguage{french}{Cadre}{Framework}}}

\newtheorem{nnconv}[thm]{\textsc{\iflanguage{french}{Convention}{Convention}}}

\newtheorem{nnterm}[thm]{\textsc{\iflanguage{french}{Terminologie}{Terminology}}}

\newtheorem{nndefn}[thm]{\textsc{\iflanguage{french}{Définition}{Definition}}}
\newenvironment{ndefn}{
    \pushQED{\qed}
    \begin{nndefn}
}{
    \popQED
    \end{nndefn}
}
\newtheorem{nnex}[thm]{\textsc{\iflanguage{french}{Exemple}{Example}}}

\newtheorem{nnexo}[thm]{\textsc{\iflanguage{french}{Exercice}{Exercise}}}

\newtheorem{nnapp}[thm]{\textsc{\iflanguage{french}{Application}{Application}}}

\newtheorem{nnrk}[thm]{\textsc{\iflanguage{french}{Remarque}{Remark}}}
\newenvironment{nrk}{
    \pushQED{\qed}
    \begin{nnrk}
}{
    \popQED
    \end{nnrk}
}
\newtheorem{nnrap}[thm]{\textsc{\iflanguage{french}{Rappel}{Reminder}}}

\newtheorem{nnnotat}[thm]{\textsc{\iflanguage{french}{Notation}{Notation}}}

\newtheorem{nnhyp}[thm]{\textsc{\iflanguage{french}{Hypothèse}{Assumption}}}

\newtheorem{nnidprf}[thm]{\textsc{\iflanguage{french}{Idée de la preuve}{Idea of the proof}}}

\theoremstyle{definition}

\incl{beamer}{\newtheorem*{abstract}{\textsc{Abstract}}}{}

\newcommand\thms{\numberwithin{thm}{section}}
\newcommand\thmss{\numberwithin{thm}{subsection}}

\newcommand\prfname{\iflanguage{french}{Preuve}{Proof}}
\makeatletter
\renewenvironment{proof}[1][\prfname]{\par\pushQED{\qed}\normalfont
  \topsep6\p@\@plus6\p@\relax
  \trivlist
  \item[\hskip\labelsep\bfseries\boldmath #1\@addpunct{.}]}{\popQED\endtrivlist\@endpefalse
}
\makeatother

  \usepackage{fancyhdr}
\usepackage{colortbl}

\arrayrulecolor{black}

\fancypagestyle{plain}{
\fancyhf{}
  \fancyhead{}
  \fancyfoot{}
  
} 
\newcommandx\vp[1][1=p]{V^{\kern-.75pt #1}}
\newcommandx\vpc[1][1=p]{\mathcal{V}^{#1}}
\newcommandx\Vp[1][1=p]{\textnormal{\textbf{V}}^{#1}} \newcommandx\Vpb[1][1=p]{\overline{\textnormal{\textbf{V}}}^{#1}} \newcommandx\vploc[1][1=p]{V^{\kern-.75pt #1}_{\textnormal{loc}}}
\newcommandx\vpcloc[1][1=p]{\mathcal{V}_{\textnormal{loc}}^{#1}}
\newcommandx\Vploc[1][1=p]{\textnormal{\textbf{V}}_{\textnormal{loc}}^{#1}} 
\newcommand\dS{d_\textnormal{\textsf{S}}}

\newcommandx\dSp[1][1=p]{d_{\textnormal{\textsf{S}},V_{#1}}}
\newcommandx\dSpo[1][1=p]{d^*_{\textnormal{\textsf{S}},V_{#1}}}
\newcommandx\dVp[1][1=p]{d_{V_{#1}}}
\newcommandx\dHa[1][1=\alpha]{d_{H_{#1}}}

\newcommandx\cvxU[1][1=]{\cvx[#1][\textnormal{\textsf{U}}]} \newcommandx\cvxS[1][1=]{\cvx[#1][\textnormal{\textsf{S}}]}
\newcommand\normi[1]{\norm{#1}_\circ}

\WithSuffix\newcommand\normi*[1]{\norm*{#1}_\circ}
\newcommand\normU[1]{\norm{#1}_{\textnormal{\textsf{U}}}}
\WithSuffix\newcommand\normU*[1]{\norm*{#1}_{\textnormal{\textsf{U}}}}
\newcommandx\normVp[2][1=p]{\norm{#2}_{V_{#1}}}
\WithSuffix\newcommandx\normVp*[2][1=p]{\norm*{#2}_{V_{#1}}}
\newcommandx\normHa[2][1=\alpha]{\norm{#2}_{\hol_{#1}}}
\WithSuffix\newcommandx\normHa*[2][1=\alpha]{\norm*{#2}_{\hol_{#1}}}

\newcommand\kl{\mathcal{K}}
\newcommand\CY{C_{\textnormal{\textsf{Y}}}}

\title{Unbounded nonconvex Young differential inclusions: existence of a measurable selection of solutions}
\author{Nathan \textsc{Benichou}}

\tikzset{bend height/.style={
        to path={.. controls ($(\tikztostart)!.25!(\tikztotarget) + (0,#1)$) 
            and ($(\tikztostart)!0.75!(\tikztotarget) + (0,#1)$) .. (\tikztotarget) \tikztonodes}
}}
\begin{document}
\pagestyle{empty}
\prettytitle
\begin{abstract}
    We study the differential inclusion  $\D z_t\in F(z_t)\D x_t$, with initial condition $z_0=\xi$,
    where $F$ is a nonconvex-valued multifunction, and $x$ a path of bounded $q$-variation,
    for some $1\lq q<2$, extending the work of Bailleul, Brault, and Coutin (2020).
    We obtain existence of local and global solutions to this inclusion without assuming $F$ bounded.
    If $z(\xi,x)$ denotes such a solution, we obtain measurability of $z$ with respect to $x$ and $\xi$.
    To establish this, we introduce a Skorokhod-type distance and prove that Young integration is continuous with respect to it.
    By the way, we prove that a compact-valued $\gamma$-Hölder map $F$ has, for any $p>1/\gamma$ and $\xi\in F(0)$,
    a selection $f(\xi)$ of bounded $p$-variation, started at $\xi$, such that $f$ is measurable in $\xi$.
\end{abstract}
\dominitoc
\printtoc
\bigskip
\smallskip

\hrule
\bigskip
\smallskip

\thmss
\pagestyle{fancy}
\setcounter{section}{-1}
\section{Introduction and notations}
\subsection{Introduction}
A differential inclusion is basically of the form
\begin{equation}
        \D z_t\in F(z_t)\D t,\quad z_0=\xi,\tag{DI}\label{di}
\end{equation}
where $F$ is a multifunction and $\xi$ an initial condition.
Differential inclusions have been extensively studied; a standard reference is \cite{aubin_differential_1984}.
Motivated by stochastic modeling, we consider the differential inclusion
\begin{equation}
        \D z_t\in F(z_t)\D x_t,\quad z_0=\xi,\tag{YDI}\label{ydi0}
\end{equation}
where $F$ is a multifunction defined on an open set $U$ of $\R^e$ with values in the nonempty compact subsets of $\R^{e\times d}$,
$x:[0,T]\to\R^d$ is a deterministic driving signal of bounded $q$-variation, for $1\lq q<2$, and $\xi\in U$ is an initial condition.
The path $x$ is typically a trajectory of some stochastic process, like a fractional Brownian motion with Hurst index $H > 1/2$.
We refer to \eqref{ydi0} as a \emph{Young differential inclusion}.
To date, two main approaches have been developed to study such inclusions.
Both are formulated under the assumption that $x$ is $\alpha$-Hölder continuous with
$1/2 < \alpha \lq 1$ and that $U = \R^e$, although they can be extended beyond this setting.
\smallskip

\textit{The approach of Bailleul, Brault and Coutin.}
In \cite{bailleul_young_2020}, a solution to \eqref{ydi0} is defined as a pair $(z,v):[0,T]\to\R^e\times\R^{e\times d}$
such that $v$ is of bounded $p$-variation for some $p\gq1$ verifying $\alpha+1/p>1$, and, for all $t\in[0,T]$,
\[z_t=\int_0^t v_\tau\D x_{\tau}\qtx{and}v_t\in F(z_t).\]
This formulation is close to the usual one for differential inclusions \cite[Ch.\,2]{aubin_differential_1984}.
Using an Euler scheme argument similar to \cite[p.\,122, Th.\,1]{aubin_differential_1984},
the authors prove existence of solutions when $F$ is $\gamma$-Hölder continuous,
for some $1/\alpha - 1 < \gamma \lq 1$, and bounded.
Apart from the boundedness assumption,
these are essentially the same conditions as in the case of Young differential equations, see \cite[Th.\,10.14]{friz_multidimensional_2010}.

\smallskip
\textit{The approach of Michta and Motyl.}
This approach applies when $F$ is convex-valued.
A solution to \eqref{ydi0} is defined as an $\alpha$-Hölder continuous path $z : [0,T] \to \R^e$ such that, for all $s \lq t$ in $[0,T]$,  
\[z_t-z_s\in\int_s^t F(z_\tau)\D x_\tau,\]
where the integral is a set-valued Young integral.
This notion plays a central role in the definition of solutions.
Several variants of such integrals have been introduced
\cite{michta_selection_2020, coutin_set-valued_2022, michta_set-valued_2022}, all inspired by Aumann's construction,
differing mainly in the class of admissible selections.

Two main results have been obtained within this framework:
\begin{itemize}
    \item In \cite{michta_solution_2023}, using the integral defined in \cite{michta_set-valued_2022},
    existence of solutions is proved under strong regularity assumptions on $F$,
    namely the existence of a locally bounded Hukuhara-derivative.
    The authors also show that the set of solutions is closed.

    \item In \cite{michta_properties_2024}, the multifunction $F$ is replaced by a set $\mathcal{H}$ of $\beta$-Hölder continuous paths from $[0,T]$ to $\R^{e \times d}$,
    with $\alpha + \beta > 1$, thereby encompassing a large class of Aumann-type integrals.
    Under suitable assumptions on $\mathcal{H}$, existence of solutions is established,
    and the solution set is shown to be bounded and closed.
    Under stronger assumptions, compactness and upper semicontinuity with respect to $x$ and $\xi$ are also obtained.
\end{itemize}

In the case of basic differential inclusions \eqref{di},
the two approaches coincide when $F$ is convex-valued and upper semicontinuous \cite[p.\,99, Lem\,1]{aubin_differential_1984}.
Whether this equivalence still holds in the Young setting remains,
to the best of our knowledge, an open question, even for convex-valued multifunctions.
\smallskip

In the present paper, we extend the results of \cite{bailleul_young_2020} in several directions:
\begin{itemize}
    \item We prove existence of global solutions in the case $U = \R^e$ under the same regularity assumptions on $F$,
    but without the boundedness condition. This yields a result analogous to the global existence theorem for Young differential equations
    \cite[Th.\,10.14]{friz_multidimensional_2010}.

    \item For a general domain $U$, we establish a local existence result,
    obtained via extension of Hölder continuous set-valued maps.
    Using similar arguments, we derive a maximal existence result when $U$ is convex.
    Removing this convexity assumption would require extending Hölder multifunctions defined on arbitrary domains,
    which remains a challenging problem.

    \item Finally, we investigate the dependence of solutions on the initial condition $\xi$ and the driving signal $x$.
    In the non-convex setting, continuity is typically out of reach; instead, we establish measurability with respect to these parameters.
    This result is particularly relevant when $x$ is a trajectory of a stochastic process
    --- we have in mind a fractional Brownian motion with Hurst index $H>1/2$ --- or when $\xi$ is a random variable.
    Although this lies beyond the scope of the present article,
    this result also opens the way to the study of a stochastic optimal control
    problem driven by a fractional Brownian motion.
\end{itemize}
The paper is organized as follows.
\begin{itemize}
    \item Section \ref{secmr} is devoted to presenting the framework of the study and stating precisely the main results.
    \item In Section \ref{sectools}, we introduce and analyze a distance $\dSp$, which combines features of the Skorokhod distance and the $p$-variation distance.
    In particular, we establish the continuity of Young integration with respect to this distance.
    This result is a key ingredient in Section \ref{secproof1},
    where we prove the measurability of solutions to \eqref{ydi} with respect to the initial condition and the driving signal.
    \item Finally, Sections \ref{secproof1} and \ref{secproof2} are devoted to the proofs of the main results stated in Section \ref{secmr}.
\end{itemize}
\subsection{Notations}\label{ssnotats}
\addssdesc{
    General notations (\pageref{not1}). Skorokhod space (\pageref{not2}). $\sigma$-algebras (\pageref{not25}).
    Subdivisions (\pageref{not3}). Hölder and $p$-variation (\pageref{not4}). Set-valued maps (\pageref{not5}). Measurable set-valued maps (\pageref{not6}).
}
\label{not1}\textit{General notations.} If $d,e\gq 1$ are integers, the set of real matrices with $e$ rows and $d$ columns is denoted by $\R^{e\times d}$.
The set of subsets of a set $X$ is denoted $\pt(X)$.
If $X$ is a topological space, we denote $\kl(X)$ the set of its compact nonempty subsets.
For a sequence $(x_n)_n\subseteq X$ we denote $\adh_n x_n$
the set of its subsequential values. When $X$ is a metric space we denote $B(x,r)$ and $\overline{B}(x,r)$
the open and closed balls of center $x\in X$ and radius $r\gq0$.
If $X$ is a vector space and $x:E\to X$ is a path defined on a subset of $\R$,
then, for $t,s\in E$ we may denote the value of $x$ at time $t$ with $x(t)$ or $x_t$,
and we write $x(s,t)=x_{s,t}=x_t-x_s$.
\bskip

\label{not2}\textit{Skorokhod space.} For $T>0$ and $X$ complete separable metric space,
the Skorokhod space of right-continuous with left limits (càdlàg) paths, denoted $\sk([0,T]\to X)$,
is equipped with its classical distance
\begin{gather*}
    \dS(x,y)=\inf_{\lambda\in\Lambda}\left(\normi{\lambda}\vee d_\infty(x,y\circ\lambda)\right),\\
    \Lambda=\{\lambda:[0,T]\to[0,T]\text{ strictly increasing and bijective}\},\\
    \normi{\lambda}=\sup_{t\ne s\in[0,T]}\abs*{\log\frac{\lambda(t)-\lambda(s)}{t-s}},
\end{gather*}
which turns it into a complete and separable metric space.
The two following propositions will be useful for us in the sequel.
\begin{npropn}\label{pskcvp}
    In $\sk([0,T]\to X)$, let $(x_n)_n$ be a sequence that converges to $x$ in distance $\dS$.
    Then, for all $t\in[0,T]$,
    $$\adh_n x_n(t)\subseteq\{x(t^-),x(t)\}.$$
\end{npropn}
\begin{npropn}\label{propsnormi}
    For every $\lambda\in\Lambda$ with $\normi{\lambda}\lq1/2$,\vspace{-.5ex}
    $$\norm{\lambda-\id}_\infty\lq 2T\normi{\lambda}\qtx{and}\lip\lambda\lq 1+2\normi{\lambda}\lq 2,$$
    where $$\lip\lambda=\sup_{t\ne s\in [0,T]}\frac{\abs{\lambda(t)-\lambda(s)}}{\abs{t-s}}.$$
\end{npropn}
For further insight on the Skorokhod space,
see \cite[\S\,12]{billingsley_convergence_2013}.
\bskip

\label{not25}\textit{$\sigma$-algebras}. Let $k\gq1$ integer and $T>0$. Unless otherwise specified, all subsets of $\R^k$ are endowed
with the Borel $\sigma$-algebra $\smash{\B(\R^k)}$. Every path space, that is, a subset of $\smash{(\R^k)^{[0,T]}}$, is endowed with the product
$\sigma$-algebra $\smash{\B(\R^k)^{\otimes[0,T]}}$, also known as cylindrical $\sigma$-algebra. It is the smallest such that every projection $\pi_t:x\mapsto x_t$, with $t\in[0,T]$,
is measurable. When $\Omega$ is some measurable space, checking that a map $f:\Omega\to(\R^k)^{[0,T]}$ is measurable
is equivalent to showing that the map $f_t\coloneqq\pi_t\circ f:\Omega\to\R^k$ is measurable for every $t\in[0,T]$.
Note that, for subsets of $\ct([0,T]\to\R^k)$, the product $\sigma$-algebra is the Borel $\sigma$-algebra of uniform distance.
More generally, for subsets of $\sk([0,T]\to\R^k)$, the product $\sigma$-algebra is the Borel $\sigma$-algebra of distance $\dS$.
\bskip

\label{not3}\textit{Subdivisions}. If $E$ is a subset of $\R$, a subdivision of $E$ is a finite subset of $E$ with at least two elements.
A subdivision denoted by $(\sigma_k)_{k\lq n}$ has $n+1$ elements, and so divides $[0,T]$
in $n$ segments. The elements of $\sigma$ are sorted in the increasing order,
so that $\sigma_0<\dots<\sigma_n$.
We can also write $[s,t]\in\sigma$ to indicate that $s,t$ are two consecutive points of $\sigma$.
The step of $\sigma$ is $\abs{\sigma}=\max_{1\lq k\lq n}(\sigma_k-\sigma_{k-1})$.
Finally we denote $\sym(E)$ the set of subdivisions of $E$, and
$\sym_a^b$ the set of subdivisions of $[a,b]$ that contains the points $a$ and $b$.
\bskip

\label{not4}\textit{Hölder, $p$-variation and oscillation}. Let $x$ be a path defined on a subset $E$ of $\R$ and valued in a normed vector space $X$.
We define the oscillation seminorm of $x$ on a subset $F\subseteq E$ by
$$\osc(x,F)=\sup_{s,t\in F}\abs{x_{s,t}}.$$
For any $0<\alpha\lq1$, we denote the $\alpha$-Hölder seminorm of $x$ on a subset $F\subseteq E$ by
$$H_\alpha(x,F)=\sup_{s\ne t\in F}\frac{\abs{x_{s,t}}}{\abs{t-s}^\alpha}.$$
For any $p\gq 1$, we denote the $p$-variation seminorm of $x$ on a subset $F\subseteq E$ by
$$V_p(x,F)=\sup_{\sigma\in\sym(F)}V_p[x,\sigma]\coloneqq\sup_{(\sigma_k)_{k\lq n}\in\sym(F)}\left(\sum_{k=1}^n\abs{x(\sigma_{k-1},\sigma_k)}^p\right)^{1/p}.$$
For more details about $p$-variation, see \cite{chistyakov_maps_1998}.
We may write $H_\alpha(x)$ and $V_p(x)$ when $F$ is the whole set $E$.
We denote respectively the space of $\alpha$-Hölder paths and the space of paths of bounded $p$-variation by
\begin{gather*}
    \cl(E\to X)=\{x:E\to X:H_\alpha(x)<\infty\},\\
    \vp(E\to X)=\{x:E\to X:V_p(x)<\infty\}.
\end{gather*}
If $q\gq p$ then $V_q(x)\lq V_p(x)$. We write
$$\vp[p+](E\to X)=\bigcap_{q>p}\vp[q](E\to X).$$
We denote by $\Vp(E\to X)$ the set of paths of bounded $p$-variation
that are càdlàg. This set is not so smaller
than $\vp(E\to X)$ because paths of bounded $p$-variation are already
regulated paths, i.e. they have left and right limits at any point.
We also denote $\vpc(E\to X)$ the set of continuous paths of bounded $p$-variation.
The space of locally $\alpha$-Hölder paths, that are $\alpha$-Hölder on every compact subset of $E$,
is denoted $\clloc(E\to X)$,
and we denote $\vploc(E\to X)$ the space of paths locally of bounded $p$-variation.
In addition to seminorms we define $p$-variation and $\alpha$-Hölder norms by
\begin{gather*}
    \normVp{x}=\norm{x}_\infty+V_p(x),\\
    \normHa{x}=\norm{x}_\infty+H_\alpha(x).
\end{gather*}
Spaces $\vp(E\to X)$ and $\cl(E\to X)$ are complete under their associated norms,
but not separable.\bskip

\label{not5}\textit{Set-valued maps}. Let $(X,\rho)$ be a metric space. For $(A_n)_n$ a sequence of subsets of $X$,
the Kuratowski limits of this sequence are
\begin{gather*}
    \Li_n A_n=\{x\in X:\exists a_n\in A_n:x=\textstyle\lim_n a_n\},\\
    \Ls_n A_n=\{x\in X:\exists a_n\in A_n:x\in\textstyle\adh_n a_n\}.
\end{gather*}
The Hausdorff distance between two nonempty subsets $A,B$ of $X$ is
$$h(A,B)=\sup_{a\in A}\rho(a,B)\vee\sup_{b\in B}\rho(b,A).$$
Let $\Omega$ be a metric space and $F:\Omega\to\pt(X)$ a set-valued map.
The domain of $F$ is the set of $\omega\in\Omega$ such that $F(\omega)\ne\emptyset$.
The preimage by $F$ of a subset $Y\subseteq X$, denoted $F^{-1}(Y)$, is the set of $\omega\in\Omega$
such that $F(\omega)\cap Y\ne\emptyset$. We say that $F$ is lower-semicontinuous
if the preimage of any open set is open,
that $F$ is upper-semicontinuous
if the preimage of any closed set is closed, and that $F$ is continuous if
it is both lower and upper semicontinuous. In that case,
for every $\omega\in\Omega$ and sequence $\omega_n\to\omega$,
we have $\overline{F}(\omega)=\Li_n F(\omega_n)$.
When nonempty-valued, $F$ is said Hausdorff-continuous if $F$ is continuous as
a single-valued map from $\Omega$ to $\pt(X)$ endowed with Hausdorff distance.
When $F$ is valued in $\kl(X)$, both notions of continuity coincide.
The notions of $\alpha$-Hölder and bounded $p$-variation
extends also to set-valued maps using Hausdorff distance.
For a more detailed exposition on continuity notions for set-valued maps, see \cite{aubin_set-valued_2009}.
\bskip

\label{not6}\textit{Measurable set-valued maps}. Let $X$ be a complete and separable metric space
and $\Omega$ a measurable space. Consider a set-valued map $F:\Omega\to\pt(X)$ with nonempty closed images.
We say that $F$ is weakly measurable if the preimage of any open set is measurable,
and that $F$ is measurable if the preimage of any closed set is measurable.
Because the space $X$ is separable, it has a countable basis of open balls,
so it suffices to check that the preimage by $F$ of any open ball of $X$ is measurable
to show that $F$ is weakly measurable.
Separability of $X$ ensures also that measurable implies weakly measurable.
If $\Omega$ is endowed with a complete $\sigma$-finite measure,
the two notions of measurability coincide, see \cite[Th.\,8.1.4]{aubin_set-valued_2009}.
Clearly a lower-semicontinuous set-valued map is weakly measurable,
and an upper-semicontinuous set-valued map is measurable.
In the context given here, every weakly measurable set-valued map $F$
has a selection $f:\Omega\to X$ that is measurable.
For a more detailed exposition on measurability notions for set-valued maps, see \cite{aubin_set-valued_2009}.

\section{Main results}\label{secmr}
\subsection{Existence of a measurable selection of solutions for unbounded nonconvex Young differential inclusions}\label{ssydi} Our framework for Young differential inclusions is essentially the one of \cite{bailleul_young_2020}.
Let $d,e\gq 1$ be two integers.
We will work with a finite time horizon $T>0$ for convenience, but it is not a restriction at all for the results
of the present article, because, by gluing solutions we can obtain the same results on ${[0,\infty[}$.
Take an open subset $U$ of $\R^e$ and a set-valued map $F:U\to\pt(\R^{e\times d})$. We denote $\Omega=U\times\vpc[q]([0,T]\to\R^d)$ for a fixed $1\lq q<2$ --- remember $\vpc[q]$
is the set of \textit{continuous} paths of bounded $q$-variation.
\begin{ndefn}[solution to a Young differential inclusion]
    Let $\omega=(\xi,x)\in\Omega$. A \textit{solution}, started from $\xi$ and defined on the interval $I={[0,T_\ast]}$, with $T_\ast\lq T$,
    to the Young differential inclusion
    \begin{equation}
        \label{ydi}\D z_t\in F(z_t)\D x_t\tag{YDI}
    \end{equation}
    is a path $z\in\vpcloc[q](I\to U)$ such that there exists
    $v\in\vploc(I\to\R^{e\times d})$, for some real number $p\gq 1$ verifying $1/p+1/q>1$,
    with, for all $t\in I$,
    \begin{gather}
        \label{ydisol} z_t=\xi+\int_0^t v_s\D x_s,\\
        v_t\in F(z_t).
    \end{gather}
    The integral in \eqref{ydisol} is a Young integral. For more details, see appendix \ref{apxa}.
\end{ndefn}
Sometimes we call solution the pair $(z,v)$, or we say that $z$ is a solution with velocity $v$.
A global solution
to \eqref{ydi} is a solution defined on the whole interval $[0,T]$.
A maximal solution is a solution which cannot be extended to a bigger interval.
\begin{nrk}
    The definition of solution we have just given could technically be extended to $q\gq 2$,
    but in this case, the path $x$ is not a well-defined rough path,
    which is an indication that things are likely to go wrong.
    For instance, existence of solutions is not guaranteed when $q\gq 2$.
\end{nrk}

The main result of the present paper is the following.

\begin{nthm}[global existence]
    \label{thydi}
    Assume $U=\R^e$ and $F$ is valued in $\kl(\R^{e\times d})$ and $\gamma$-Hölder
    for some $q-1<\gamma\lq 1$. Then, for every $\omega=(\xi,x)\in\Omega$, there exists
    a global solution $\smash{z(\omega)\in\vpc[q]([0,T]\to U)}$
    to \eqref{ydi} started at $\xi$ with velocity $\smash{v(\omega)\in\Vp[q/\gamma+]([0,T]\to\R^{e\times d})}$.
    Moreover, we can choose the latter in a way that the maps $z$ and $v$ are measurable
    from $\Omega$ to their target path spaces.
\end{nthm}

Theorem \ref{thydi} differs from \cite[Th.\,1.1]{bailleul_young_2020} in two ways.
First, the map $F$ is no longer assumed to be bounded.
Second, the maps $z$ and $v$ are measurable with respect
to $x$ and $\xi$.
\smallskip

Complementing the preceding global existence result,
the next two theorems establish local existence.
\begin{nthm}[local existence]
    \label{thydiloc}
    Assume that $F$ is valued in $\kl(\R^{e\times d})$ and locally $\gamma$-Hölder
    for some $q-1<\gamma\lq 1$. Then, for all $\omega=(\xi,x)\in\Omega$, there exists
    $T_\ast(\omega)>0$
    and a local solution $\smash{z(\omega)\in\vpc[q]([0,T_\ast(\omega)]\to U)}$
    to \eqref{ydi} started at $\xi$ with velocity $\smash{v(\omega)\in\Vp[q/\gamma+]([0,T_\ast(\omega)]\to\R^{e\times d})}$.
    Moreover, we can choose the latter in a way that the maps $T_\ast$ and $(z,v)\ind_{[0,T_\ast]}$ are measurable.
\end{nthm}
\begin{nthm}[maximal existence]
    \label{thydimax}
    Suppose $U$ is convex and assume that $F$ is valued in $\kl(\R^{e\times d})$ and locally $\gamma$-Hölder
    for some $q-1<\gamma\lq 1$. Then, for every $\omega=(\xi,x)\in\Omega$, there exists
    an interval $I(\omega)={[0,T_\ast(\omega)]}$, with $T_\ast(\omega)>0$,
    and a local solution $\smash{z(\omega)\in\vpcloc[q](I(\omega)\to U)}$
    to \eqref{ydi}, started at $\xi$ and with velocity $\smash{v(\omega)\in\Vploc[q/\gamma+](I(\omega)\to\R^{e\times d})}$,
    with the property that $z_{I(\omega)}(\omega)$ is not precompact in $U$ if $I(\omega)\ne [0,T]$.
    In particular the latter is a maximal solution.
    Moreover, we can choose it in a way that the maps $T_\ast$ and $(z,v)\ind_{I}$ are measurable.
\end{nthm}
\subsection{Selection of an Hölderian set-valued map with measurability with respect to initial value}
While writing the article \cite{bailleul_young_2020}, the authors realized that their
proof of existence of solutions to \eqref{ydi} can easily be adapted in order to prove
that any $\gamma$-Hölder set-valued map $F:[0,T]\to\kl(\R^d)$
has a selection in $\Vp[1/\gamma+]([0,T]\to\R^d)$ started at an arbitrary point $\omega\in F(0)$.
This have been a partial answer to the open question \cite[Rk.\,8.1]{chistyakov_maps_1998}.
In turn, we can also adapt the proof of \ref{thydi}, resulting in the following theorem,
which is the same selection result with some additional measurability.
\begin{nthm}\label{thsel}
    Let $F:[0,T]\to\kl(\R^d)$ be $\gamma$-Hölder for some $0<\gamma\lq1$. Denote $\Omega=F(0)$. For all $\omega\in\Omega$, there exists a selection $f(\omega)\in\Vp[1/\gamma+]([0,T]\to\R^d)$ of $F$
    such that $f_0(\omega)=\omega$.
    We also have the
    estimate, for all $p>1/\gamma$,
    \begin{equation}\label{estselvpor}
        V_p(f(\omega))\lq C H_\gamma(F), \end{equation}
    with $C$ a constant depending on $p$ and $\gamma$. Moreover, we can choose $f(\omega)$ such that $f$ is measurable from $\Omega$ to its target path space.
\end{nthm}
Theorem \ref{thydi} is proved in Section \ref{secproof1}.
The other proofs are given in \ref{secproof2}.
\section{Continuity of Young integration under a Skorokhod-type distance}\label{sectools}
\thms
Let $p\gq1$ and $0<\alpha\lq1$ such that $\alpha+1/p>1$.
Fix integers $d,e\gq1$, a finite time horizon $T>0$ and a path $x\in\cl([0,T]\to\R^d)$.
We know from \ref{pyibil} that the map
\begin{equation}
    \Vp([0,T]\to \R^{e\times d})\ni v\longmapsto \int_0^\cdot v\D x\in\cl([0,T]\to \R^e)\label{indefyi}
\end{equation}
is continuous under norms $\normd_{V_p}$ and $\normd_{H_\alpha}$.
The goal of this section is to show that this continuity holds when replacing the norm
$\normd_{V_p}$ by a weaker distance of Skorokhod type --- at the cost of a little weakening of the norm $\normd_{H_\alpha}$.
This will be a tool we will use in \ref{sspglob} for proving Theorem \ref{thydi}.
The genuine Skorokhod distance $\dS$ is not adapted to our problem. We
shall define a stronger distance
$$\dSp(x,y)=\inf_{\lambda\in\Lambda}\left(\normi{\lambda}\vee \normVp{x-y\circ\lambda}\right),$$
for any two paths $x,y:[0,T]\to X$, where $X$ is a separable Banach space.
The norm $\normd_{V_p}$ is invariant under nondecreasing bijective changes of time,
hence $\dSp$ is indeed a distance on $\Vp([0,T]\to X)$, and
$$\dS(x,y)\lq \dSp(x,y)\lq \normVp{x-y}.$$
The distance $\dSp$ has been used in \cite[\S\,5.1]{friz_differential_2018}
to derive some continuity results on rough differential equations. Many of the properties of the space $(\sk,\dS)$
do not seem to carry over to $(\Vp,\dSp)$, or at least not in any straightforward way.
For instance, we do not know whether $\Vp$ is separable under $\dSp$.
The situation improves slightly when $\Vp$ is regarded as a subset of $\Vp[q]$
for some $q>p$, a fact that can be attributed primarily to interpolation inequalities. Let us recap in a few lines what interpolation inequalities are.
For any path $x:[0,T]\to X$, for $q>p$ and $0<\beta<\alpha$, one can check that
\begin{align}V_q(x)     & \lq\osc(x)^{1-{p/q}}V_p(x)^{p/q},                        \\
             H_\beta(x) & \lq\osc(x)^{1-{\beta/\alpha}}H_\alpha(x)^{\beta/\alpha}.\end{align} Then if $y:[0,T]\to X$ is another path, we deduce
\begin{align}V_q(x-y)                           & \lq(2\norm{x-y}_\infty)^{1-{p/q}}(V_p(x)+V_p(y))^{p/q},                             \\
             \label{ineqinterpham1}H_\beta(x-y) & \lq(2\norm{x-y}_\infty)^{1-{\beta/\alpha}}(H_\alpha(x)+H_\alpha(y))^{\beta/\alpha}.\end{align}
And these inequalities extend naturally to norms:
\begin{align}\label{ineqinterpvp}\normVp[q]{x-y}     & \lq\norm{x-y}_\infty^{1-{p/q}}\big(\norm{x-y}_\infty^{p/q}+2V_p(x)^{p/q}+2V_p(y)^{p/q}\big),                                               \\
             \label{ineqinterpha}\normHa[\beta]{x-y} & \lq\norm{x-y}_\infty^{1-{\beta/\alpha}}\big(\norm{x-y}_\infty^{\beta/\alpha}+2H_\alpha(x)^{\beta/\alpha}+2H_\alpha(y)^{\beta/\alpha}\big).\end{align}
There is a Skorokhod counterpart of inequality \eqref{ineqinterpvp},
it is given in the following lemma.
Its proof is a direct computation.
\begin{nlem}\label{pinterpineqsk}
    For any two paths $x,y:[0,T]\to X$ and all $q>p$,
    \begin{equation}
        \dSp[q](x,y)\lq\dS(x,y)^{1-{p/q}}\left(\dS(x,y)^{p/q}+2V_p(x)^{p/q}+2V_p(y)^{p/q}\right).
    \end{equation}
\end{nlem}
Now we can show the main result of this section.
\begin{nthm}\label{contiy}
    The map given in \eqref{indefyi} is continuous under distance $\dSp$ and norm $\normd_{H_\beta}$\,
    for any $0<\beta<\alpha$.
\end{nthm}
\begin{proof}
    Let $v,w\in \Vp([0,T]\to \R^{e\times d})$ be two paths such that $\dSp(v,w)<1/2$.
    Take $\epsilon>0$ and a change of time $\lambda\in\Lambda$ verifying
    $$\normi{\lambda}\vee\normVp{v\circ\lambda-w}<\dSp(v,w)+\epsilon<\tfrac{1}{2}.$$
    For $t\in[0,T]$,
    $$A(t)\coloneqq\int_0^t v\D x-\int_0^t w\D x=f(t)+g(t)+h(t),$$
    where
    \begin{align*}
        f(t) & =\int_0^t v\D x-\int_0^t (v\circ\lambda)\D (x\circ\lambda),               \\
        g(t) & =\int_0^t (v\circ\lambda)\D (x\circ\lambda)-\int_0^t (v\circ\lambda)\D x, \\
        h(t) & =\int_0^t (v\circ\lambda)\D x-\int_0^t w\D x.
    \end{align*}
We will bound $\norm{f}_\infty$, $\norm{g}_\infty$ and $H_\alpha(h)$ in order to bound $H_\beta(A)$ and finally $\normHa[\beta]{A}$.
    By \ref{pyicv}, \ref{pyivp} and \ref{propsnormi}, for all $t\in[0,T]$,
    \begin{align*}
        \abs{f(t)}=\abs*{\int_0^t v\D x-\int_0^{\lambda(t)}v\D x}&=\abs*{\int_{\lambda(t)}^t v\D x}\\
         & \lq \CY(1/\alpha,p)\normVp{v}H_\alpha(x)\abs{t-\lambda(t)}^\alpha        \\
                                                                 & \lq 2^\alpha \CY(1/\alpha,p)H_\alpha(x)T^\alpha\normVp{v}\normi{\lambda}^\alpha.
    \end{align*}
    By \ref{pyibil} and \ref{pyivp}, for all $t\in[0,T]$,
    $$\abs{g(t)}=\abs*{\int_0^t (v\circ\lambda)\D (x\circ\lambda-x)}\lq \CY(1/\gamma,p)\normVp{v}H_\gamma(x\circ\lambda-x)T^\gamma,$$
    for any $0<\gamma<\alpha$ such that $\gamma+1/p>1$. By inequality \eqref{ineqinterpham1} and by \ref{propsnormi},
    \begin{align*}
        H_\gamma(x\circ\lambda-x) & \lq (2\norm{x\circ\lambda-x}_\infty)^{1-{\gamma/\alpha}}(H_\alpha(x\circ\lambda)+H_\alpha(x))^{\gamma/\alpha}               \\
                                  & \lq (2H_\alpha(x)\norm{\lambda-\id}_\infty^\alpha)^{1-{\gamma/\alpha}}(H_\alpha(x)((\lip\lambda)^\alpha+1))^{\gamma/\alpha} \\
                                  & \lq 12H_\alpha(x)T^{\alpha-\gamma}\normi{\lambda}^{\alpha-\gamma}.
    \end{align*}
    We deduce
    $$\abs{g(t)}\lq 12\CY(1/\gamma,p)H_\alpha(x)T^\alpha\normVp{v}\normi{\lambda}^{\alpha-\gamma}.$$
    By \ref{pyibil} and \ref{pyiha},
    \begin{align*}
        H_\alpha(h)=H_\alpha\left(\int_0^\cdot (v\circ\lambda-w)\D x\right) & \lq \CY(1/\alpha,p)H_\alpha(x)\normVp{v\circ\lambda-w}.
    \end{align*}
    Now, let $0<\beta<\alpha$. Using \eqref{ineqinterpham1}, we have by inequalities above,
    \begin{align*}
        H_\beta(f+g) & \lq(2\norm{f+g}_\infty)^{1-{\beta/\alpha}}(H_\alpha(f)+H_\alpha(g))^{\beta/\alpha}                                             \\
                     & \lq C_1(H_\alpha(x)T^\alpha\normVp{v}\normi{\lambda}^{\alpha-\gamma})^{1-{\beta/\alpha}}(\normVp{v}H_\alpha(x))^{\beta/\alpha} \\
                     & \lq C_1H_\alpha(x)T^{\alpha-\beta}\normVp{v}\normi{\lambda}^{(\alpha-\gamma)(1-\beta/\alpha)},
    \end{align*}
    for some constant $C_1$ depending on $\alpha$, $\gamma$, $\beta$ and $p$.
    By basic properties of Holder seminorms we also have
    $$H_\beta(h)\lq T^{\alpha-\beta}H_\alpha(x)\lq \CY(1/\alpha,p)H_\alpha(x)T^{\alpha-\beta}\normVp{v\circ\lambda-w}.$$
    Summing up, there is some constant $C_2$ depending on $\alpha$, $\gamma$, $\beta$ and $p$, such that
    \begin{align*}
        H_\beta(A)
         & \lq C_2H_\alpha(x)T^{\alpha-\beta}(\normVp{v}\vee1)(\normi{\lambda}\vee\normVp{v\circ\lambda-w})^{(\alpha-\gamma)(1-\beta/\alpha)} \\
         & \lq C_2H_\alpha(x)T^{\alpha-\beta}(\normVp{v}\vee1)(\dSp(v,w)+\epsilon)^{(\alpha-\gamma)(1-\beta/\alpha)}.
    \end{align*}
    Letting $\epsilon\to0$ we obtain$$H_\beta(A)\lq C_2H_\alpha(x)T^{\alpha-\beta}(\normVp{v}\vee1)\dSp(v,w)^{(\alpha-\gamma)(1-\beta/\alpha)},$$
    from which we deduce the estimate
    \begin{equation}
        \normHa[\beta]{A}\lq C_2H_\alpha(x)T^{\alpha-\beta}(T^\beta+1)(\normVp{v}\vee1)\dSp(v,w)^{(\alpha-\gamma)(1-\beta/\alpha)}.\label{majyoungdsp}
    \end{equation}
    Considering $v$ fixed, the expression \eqref{majyoungdsp} tends to $0$ as $\dSp(v,w)\to0$. Hence the desired continuity.
\end{proof}
\section{Proof of the global existence result}\label{secproof1}
\thms
This section is dedicated to proving Theorem \ref{thydi}.
As a first step we do in Subsection \ref{sspglob}
a proof in the restricted case of an $\alpha$-Hölder
$x$ verifying $H_\alpha(x)\lq1$ for some $1/2<\alpha\lq 1$.
This latter proof makes use of estimates of Section \ref{sspest}.
The general proof is then obtained in Subsection \ref{sspglobgen}
by splitting a path $x\in\vpc[q]([0,T]\to\R^d)$
into a composition of a $1/q$-Hölder path and a nondecreasing path, like in \cite[\S\,3]{chistyakov_maps_1998}.
\bskip

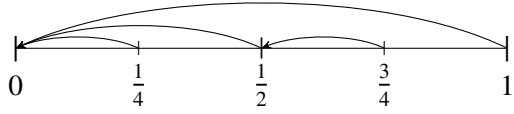
\begin{figure}[h]
    \begin{center}
        \begin{tikzpicture}[x=.5\textwidth]
            \draw (0,0) -- (1,0);
            \draw (0,0) node {$|$} node[yshift=-1.25em] {$0$}
            (1/4,0) node {\tiny$|$} node[yshift=-1.25em] {$\frac{1}{4}$}
            (1/2,0) node {\footnotesize$|$} node[yshift=-1.25em] {$\frac{1}{2}$}
            (3/4,0) node {\tiny$|$} node[yshift=-1.25em] {$\frac{3}{4}$}
            (1,0) node {$|$} node[yshift=-1.25em] {$1$};
            \path[->,>=stealth'] (1,0) edge[bend height=.8] (0,0);
            \path[->,>=stealth'] (1/2,0) edge[bend height=.4] (0,0);
            \path[->,>=stealth'] (1/4,0) edge[bend height=.2] (0,0);
            \path[->,>=stealth'] (3/4,0) edge[bend height=.2] (1/2,0);
        \end{tikzpicture}
    \end{center}
    \vspace{-1em}
    \small\caption{Ancestors of the dyadic numbers of\\order less than $2$ of the interval $[0,1]$.}
    \label{figances}
\end{figure}

We will make extensive use of dyadic subdivisions
so let us give the notations. For $S>0$ and every
$m\in\N$, the $m^\textnormal{th}$ dyadic subdivision of the interval $[0,S]$ is the set $$D^m(S)=\{kS2^{-m}\,|\,k=0,\dots,2^m\}.$$
It is a subdivision of the interval $[0,S]$ in $2^m$ intervals as defined in Subsection \ref{ssnotats}.
According to the notations we gave for subdivisions, we denote $D^m_k(S)=kS2^{-m}$ for each $k=0,\dots,2^m$.
The increasing union $D(S)=\bigcup_m D^m(S)$ is called the set of dyadic numbers of the interval $[0,S]$.
For every $t\in D(S)$, the order of $t$, denoted by $M(t)$, is the smallest integer $m\in\N$
such that $t\in D^m(S)$. Note that for every $m\in\N$,
$D^m(S)$ is the set of dyadic numbers of order less than $m$.
For $t\in D(S)\ssm\{0\}$, the ancestor of $t$ is
\begin{align*}
    a(t) & =\max\{s\in D^{M(t)}(S):s<t\} \\&=\max\{s\in D^{M(t)-1}(S):s<t\}=t-S2^{-M(t)}.
\end{align*}
The number $a(t)$ is in $D^{M(t)-1}$ but not always of order $M(t)-1$. For instance,
the number $a(D^2_1(S))=0$ is of order $0$ but $D^2_1(S)$ is of order $2$.
Finally, if $t$ is any element of $D^m(S)$ for $m\in\N$, then $s=t+S2^{-m-1}$ is of order $m+1$ and $a(s)=t$.
See figure \ref{figances} for an illustration of the notion of ancestor.

Before starting the proofs we exhibit a useful Ascoli-type theorem \cite[p.\,15]{aubin_differential_1984}.
\begin{ndefn}[equioscillating familly]
    Let $\mathcal{F}$ be a subset of the set of bounded functions from a set $E$ to
    a metric space $X$. We say that $\mathcal{F}$ is \textit{equioscillating} if, for every $\epsilon>0$,
    there exists a partition $(E_i)_{i\lq n}$ of $E$ such that $\osc(x,E_i)\lq\epsilon$ for every $x\in\mathcal{F}$ and all $i\lq n$.
\end{ndefn}
\begin{cthm}{ASC}[discontinuous Arzelà-Ascoli]\label{thasc}
    Let $\mathcal{F}$ a subset of the set of bounded functions from a set $E$ to
    a metric space $X$ verifying
    \begin{enumtimes}
        \item $\mathcal{F}$ is equioscillating.
        \item For all $t\in E$, the set $\{x(t)\,|\,x\in\mathcal{F}\}\subseteq X$ is totally bounded.
    \end{enumtimes}
    Then, under uniform distance, $\mathcal{F}$ is totally bounded.
    In particular, if $X$ is complete, $\mathcal{F}$ is precompact in uniform distance.
\end{cthm}
\thmss
\subsection{Restricted case}\label{sspglob}
Fix some $1/2<\alpha\lq1$ and denote
$\tstc=\{x:[0,T]\to\R^d:H_\alpha(x)\lq 1\}$.
In this subsection we prove the following.
\begin{nlem}\label{lemrest}
    Theorem \ref{thydi} is true in the restricted case of \,$\Omega=U\times \tstc$.
\end{nlem}
So let $U=\R^e$ and take a set-valued map $F:U\to \kl(\R^{e\times d})$ that is $\gamma$-Hölder
for some $1/\alpha-1<\gamma\lq 1$. Define $\Omega=U\times \tstc$.

\bskip

For a subset $A$ of a normed vector space we denote $\norm{A}=\sup_{a\in A}\abs{a}$. By definition of the Hausdorff distance, we have, for any subsets $A$ and $B$,
$$\norm{B}\lq\norm{A}+h(A,B).$$
In particular, the map $U\ni\xi\mapsto\norm{F(\xi)}$ is $\gamma$-Hölder
of constant less than $H_\gamma(F)$.
\bskip

The proof we present is a refinement of the proof of \cite[Th.\,1.1]{bailleul_young_2020}
in which the constants are tracked more carefully and measurability issues are handled explicitly at each step.
It is a proof by Euler scheme,
which is somewhat classical in differential inclusions. One can find a prototype
in \cite[p.\,112]{aubin_differential_1984}.
The first step, which is the more technical, consists in constructing
many approximate solutions. It is the content of the two next propositions.
To prove these propositions will make use of the results of Lemma \ref{techlem1},
which we have placed in a separate subsection for clarity.
\bskip

In all this subsection, we fix some real numbers $1/(1+\gamma)<\beta<\alpha$, $0<\theta<1$ and $1/(\gamma\alpha)<p_0<1/(1-\alpha)$.
We then set constants
\begin{align*}
    C_0 & =T\wedge2^{-\frac{1}{\theta(\alpha-\beta)}}\wedge\left(\frac{H_\gamma(F)}{1-2^{1-(\gamma+1)\alpha}}\right)^{-\frac{1}{\theta(\alpha-\beta)}}
    \wedge\left(\frac{2H_\gamma(F)}{1-2^{-\alpha\gamma}}\right)^{-\frac{1}{\gamma(1-\theta)(\alpha-\beta)}},                                            \\[1ex]
    C_1 & =\CY(1/\alpha,p_0)\left(1+C_2(p_0)\right),\\
    C_2(p) & =\left(\frac{2}{1-2^{1-p\gamma\alpha}}\right)^{1/p},
\end{align*}
where $\CY$ is defined in \ref{pyivp}.
\begin{npropn}[existence of local approximate solutions]\label{pexistloc} Consider a map $R:\Omega\to{]0,\infty[}$ and denote, for all $\omega=(\xi,x)\in\Omega$,
    \begin{gather*}
        K(\omega)=1+\norm{F(\xi)}+H_\gamma(F)R(\omega)^\gamma,\label{defK}\\
        T_0(\omega)=C_0\wedge\left(\frac{R(\omega)}{C_1 K(\omega)}\right)^{1/\alpha}.\label{defT}
    \end{gather*}
    If $T_\ast:\Omega\to{]0,\infty[}$ is any map such that $T_\ast\lq T_0$,
    then, for every $m\in\N$ and all $\omega=(\xi,x)\in\Omega$, there exist two paths $v^m(\omega)\in\Vp[1/(\gamma\alpha)+]([0,T_\ast(\omega)]\to\R^{e\times d})$
    and $z^m(\omega)\in \cl([0,T_\ast(\omega)]\to U)$ such that
    \begin{enumtimesrom}
        \item\label{locca} $z^m_t(\omega)=\xi+\int_0^t v^m_s(\omega)\D x_s$ \ for all $t\in[0,T_\ast(\omega)]$.
        \item\label{loccb} For every $t\in[0,T_\ast(\omega)]$, we have $\abs{z^m_t(\omega)-\xi}\lq R(\omega)$ and $\abs{v^m_t(\omega)}\lq K(\omega)$.
        \item\label{loccc} $v^m_t(\omega)\in F(z^m_t(\omega))$ \ for every\footnotemark\ $t\in D^m(T_\ast(\omega))$.
        \footnotetext{$D^m$ is defined at the begining of the current section.}
        \item\label{loccd} We have the estimates, for $n\in\N$ and $p>1/(\gamma\alpha)$,
        \begin{gather}
            \sup_{[s,t]\in D^n(T_\ast(\omega))}\osc\left(v^m(\omega),{[s,t[}\right)\lq2^{-\gamma\alpha n}K(\omega)^\gamma T_\ast(\omega)^{\gamma\beta}\ind_{n<m},\label{realestosc}\\ V_p(v^m(\omega),{[0,T_\ast(\omega)[})\lq C_2(p)K(\omega)^\gamma T_\ast(\omega)^{\gamma\beta}.\label{realestvp}
        \end{gather}
        \item\label{locce} Every map
        $$\Omega\ni\omega\longmapsto \left(z^m_{tT_\ast(\omega)}(\omega),v^m_{tT_\ast(\omega)}(\omega)\right),$$
        with $t\in[0,1]$, is measurable.
    \end{enumtimesrom}
\end{npropn}
\begin{proof}
    In the whole proof, $D^m(\omega)$ will stand for $D^m(T_\ast(\omega))$. Fix an integer $m\in\N$.
The values of $F$ are closed and nonempty subsets of the Polish space $\R^{e\times d}$.
    Moreover, $F$ is $\gamma$-Hölder and therefore is measurable,
    so it has a measurable selection, that we denote by $\phi$.
    We then set, for all $\omega=(\xi,x)\in\Omega$,
    \begin{align*}
        v^m_t(\omega) & =\phi(\xi)               & \hspace{-5em}\text{for all $t\in[0,D^m_1(\omega)[$}, \\
        z^m_t(\omega) & =\xi+\phi(\xi)x_{0,t} & \hspace{-5em}\text{for all $t\in[0,D^m_1(\omega)]$}.
    \end{align*}
    The first step of the construction of $v^m(\omega)$ and $z^m(\omega)$ is done.
    We now perform the next ones by induction on $k$. Suppose that for all $\omega=(\xi,x)\in\Omega$ we have constructed the paths $v^m(\omega)$ and $z^m(\omega)$
    respectively on the intervals ${[0,\tau(\omega)[}$ and $[0,\tau(\omega)]$, with $\tau(\omega)\coloneqq D^m_k(\omega)$
    for some $1\lq k< 2^m$, and that
    \begin{enumtimes}
        \item\label{loccpa} $z^m_t(\omega)=\xi+\int_0^t v^m_s(\omega)\D x_s$ \ for all $t\in[0,\tau(\omega)]$.
        \item\label{loccpb} $v^m_t(\omega)\in F(z^m_t(\omega))$ \ for every $t\in D^m(\omega)\cap{[0,\tau(\omega)[}$.
        \item\label{loccpc} $v^m(\omega)$ is constant on every interval $[D^m_{i-1}(\omega),D^m_i(\omega)[$ with $i=1,\dots,k$.
                        \item\label{loccpd} $v^m(\omega)$ is bounded by $K(\omega)$.
                        \item\label{loccpe} For all $t\in D^m(\omega)\cap{]0,\tau(\omega)[}$, we have $$\abs{v^m_{t}(\omega)-v^m_{a(t)}(\omega)}\lq H_\gamma(F)\abs{z^m_{t}(\omega)-z^m_{a(t)}(\omega)}^\gamma.$$
        \item\label{loccpf} Every map
        $$\Omega\ni\omega\mapsto v^m_{tT_\ast(\omega)}(\omega)
            \qtx{and}\Omega\ni\omega\mapsto z^m_{tT_\ast(\omega)}(\omega),$$ with respectively $t\in[0,D^m_k(1)[$ and $t\in[0,D^m_k(1)]$, is measurable.
    \end{enumtimes}
    For $\omega=(\xi,x)\in\Omega$, we can extend the path $v^m(\omega)$ into a piecewise constant path \mbox{$\tilde{v}:[0,T_\ast(\omega)]\to\R^{e\times d}$}
    setting
    $$\tilde{v}_s(\omega)=v_s(\omega) \text{ for every $s\in D^m(\omega)\cap {[0,\tau(\omega)[}$},$$
    and then recursively
    $$\tilde{v}_s(\omega)=\tilde{v}_{a(s)}(\omega)\text{ for every $s\in D^m(\omega)\cap [\tau(\omega),T_\ast(\omega)]$}.$$
    Applying Lemma \ref{techlem1} on $\tilde{v}$ we obtain the estimate
    $$V_{p_0}(v^m(\omega),{[0,\tau(\omega)[})\lq C_2(p_0) K(\omega)^\gamma T_\ast(\omega)^{\gamma\beta}.$$
    Thus by \ref{pyivp} and \eqref{defT}, for all $t\in [0,\tau(\omega)]$, \begin{equation}\label{zr}
        \begin{aligned}
            \abs{z^m_{0,t}(\omega)} & \lq \CY(1/\alpha,p_0)(K(\omega)+C_2(p_0)K(\omega)^\gamma T_\ast(\omega)^{\gamma\beta})T_\ast(\omega)^\alpha \\
                                     & \lq C_1K(\omega)T_\ast(\omega)^\alpha\\
                                     & \lq R(\omega).
        \end{aligned}
    \end{equation}
    Let us now construct the next values of $v^m(\omega)$ and $z^m(\omega)$. As $F$ is $\gamma$-Hölder we have
    $$h\big(F\big(z^m_{\tau(\omega)}(\omega)\big),F\big(z^m_{a(\tau(\omega))}(\omega)\big)\big)\lq H_\gamma(F)\big|z^m_{\tau(\omega)}(\omega)-z^m_{a(\tau(\omega))}(\omega)\big|^\gamma\eqcolon r(\omega).$$
    On the other hand, we have $v^m_{a(\tau(\omega))}(\omega)\in F(z^m_{a(\tau(\omega))}(\omega))$ by \ref{loccpb}. Then, because $F$ is closed-valued,
    the set
    $$G(\omega)\coloneqq\overline{B}\big(v^m_{a(\tau(\omega))}(\omega),r(\omega)\big)\cap F\big(z^m_{\tau(\omega)}(\omega)\big)$$
    is nonempty and closed. To show that $G$ is measurable, we write, for a closed subset $A$ of $\R^{e\times d}$,
    $G^{-1}(A)=g^{-1}(H^{-1}(A))$, where
    \begin{gather*}
        g:\Omega\ni\omega\longmapsto\left(z^m_{\tau(\omega)}(\omega),v^m_{a(\tau(\omega))}(\omega),r(\omega)\right),\\
        H:\Omega\times\R^{e\times d}\times{[0,\infty[}\ni (z,v,r)\longmapsto \overline{B}(v,r)\cap F(z).
\end{gather*}
    By \ref{loccpf} the map $g$ is measurable, and by \cite[p.\,41]{aubin_differential_1984} the set-valued map $H$ is upper-semicontinuous on its domain,
    hence $H^{-1}(A)$ is closed and $G^{-1}(A)$ is measurable.
    We are therefore able to take a measurable selection of $G$ which we denote by $\phi$.
    We then set, for all $\omega=(\xi,x)\in\Omega$,
    \begin{align*}
        v^m_t(\omega) & =\phi(\omega)                                              & \hspace{-1.5em}\text{for all $t\in[D^m_k(\omega),D^m_{k+1}(\omega)[$}, \\
        z^m_t(\omega) & =z^m_{\tau(\omega)}(\omega)+\phi(\omega)x_{\tau(\omega),t} & \hspace{-1.5em}\text{for all $t\in[D^m_k(\omega),D^m_{k+1}(\omega)]$}.
    \end{align*}
    Propositions \ref{loccpa}, \ref{loccpb}, \ref{loccpc}, \ref{loccpe} and \ref{loccpf} are clearly true at rank $k+1$.
    To prove \ref{loccpe}, we use the definition of $\phi$ and \eqref{zr} to obtain
    \begin{align*}
        \abs{\phi(\omega)} & \lq\norm*{F(z^m_{\tau(\omega)}(\omega))}                                  \\
                           & \lq\norm{F(\xi)}+H_\gamma(F)\abs{z^m_{0,\tau(\omega)}(\omega)}^\gamma \\
                           & \lq\norm{F(\xi)}+H_\gamma(F)R(\omega)^\gamma                           \\
                           & \lq K(\omega).
    \end{align*}
    By running this induction, we build paths $v^m(\omega)$ and $z^m(\omega)$ respectively on ${[0,T_\ast(\omega)[}$ and $[0,T_\ast(\omega)]$.
    To finish, we set $v^m_{T_\ast(\omega)}(\omega)$ as a measurable selection of $F(z^m_{T_\ast(\omega)}(\omega))$.
    By construction, propositions \ref{locca}, \ref{loccc} and \ref{locce}
    are true. $v^m(\omega)$ is definitely bounded by $K(\omega)$ and, by Lemma \ref{techlem1},  we have \ref{loccd}.
    Finally by the same computation as in \eqref{zr} we obtain \ref{loccb} and conclude
    the proof of Proposition \ref{pexistloc}.
\end{proof}
\begin{npropn}[existence of global approximate solutions]\label{pexistglob}
    For every $m\in\N$ and all $\omega=(\xi,x)\in\Omega$, there exist two paths $v^m(\omega)\in\Vp[1/(\gamma\alpha)+]([0,T]\to\R^{e\times d})$
    and $z^m(\omega)\in \cl([0,T]\to U)$, such that
    \begin{enumtimesrom}
        \item\label{globca} $z^m_t(\omega)=\xi+\int_0^t v^m_s(\omega)\D x_s$ \ for all $t\in[0,T]$.
        \item\label{globcb} There exists a nondecreasing familly $(\T_m)_m$ of subsets of $[0,T]$
        containing $0$ and $T$, not depending on $\omega$, such that $v^m_t(\omega)\in F(z^m_t(\omega))$ for all $t\in\T_m$,
        and the set $\T=\bigcup_m\T_m$ is dense in $[0,T]$.
        \item\label{globcc} The familly $(v^m(\omega))_m$ is equioscillating.
        \item\label{globce} $V_p(v^m(\omega))$ is bounded uniformly in $m$ for every $p>1/(\gamma\alpha)$.
        \item\label{globcd} The maps $z^m$ and $v^m$ are measurable from $\Omega$ to their target path spaces.
    \end{enumtimesrom}
\end{npropn}
\begin{proof}
    This proof consists in constructing a global solution by gluing local ones
    that are defined on a time span independent of $\omega$.
    An illustration of the quantities involved is given on figure \ref{figglob}.
    For $\omega=(\xi,x)\in\Omega$
    we denote
    $$R(\omega)=R(\xi)=1\vee\left(\frac{1+\norm{F(\xi)}}{H_\gamma(F)}\right)^{1/\gamma},$$
    so that $R$ is positive, continuous in $\xi$, and
    $$\frac{R(\xi)}{C_1(1+\norm{F(\xi)}+H_\gamma(F)R(\xi)^\gamma)}\gq\frac{1}{2C_1H_\gamma(F)}.$$
Next, we define maps $K$, $T_0$ as in Proposition \ref{pexistloc}.
    Note that $K$ and $T_0$ depend only on $\xi$.
    We have chosen $R(\xi)$ in a way that $T_0(\xi)\gq C_3$ for every $\xi\in U$, where $C_3=C_0\wedge(2C_1H_\gamma(F))^{-\alpha}$.
    We thus choose an integer $k\gq1$ large enough so that $\tilde{T}\coloneqq T/k\lq C_3$,
    and apply Proposition \ref{pexistloc} with $T_\ast\equiv \tilde{T}$. This way we obtain two sequences of maps $\tilde{v}$ and $\tilde{z}$.
    We define now, for every $m\in\N$, $\omega=(\xi,x)\in\Omega$, $1\lq i\lq k$ and $t\in{[\smash{(i-1)\tilde{T},i\tilde{T}}[}$,
    closing this interval if $i=k$,
    \begin{align*}
        z^m_t(\omega) & =\tilde{z}^m_s(\cdot,x)\circ\big(\tilde{z}^m_{\tilde{T}}(\cdot,x)\big)^{\circ \,i-1}(\xi), \\
        v^m_t(\omega) & =\tilde{v}^m_s(\cdot,x)\circ\big(\tilde{v}^m_{\tilde{T}}(\cdot,x)\big)^{\circ \,i-1}(\xi),
    \end{align*}
    where $s=t-(i-1)\tilde{T}$. Using \ref{pexistloc} we see clearly that propositions \ref{globca} and \ref{globcd} are true, and proposition \ref{globcb}
    is verified with $\T_m=\{i\tilde{T}2^{-m}\,|\,i=0,\dots,k2^m\}$. To prove proposition \ref{globcc},
    let $n\in\N$. For any two consecutive points $s,t$ of $\T_n$
    and any $m\in\N$, we have, by \eqref{realestosc},
    $$\osc\left(v^m(\omega),{[s,t[}\right)\lq2^{-\gamma\alpha n}\overline{K}^m(\omega)^\gamma \tilde{T}^{\gamma\beta}$$
    with
    $$\overline{K}^m(\omega)=\max_{i=1}^k K(\xi^m_i)\qtx{and}\xi^m_i=z^m_{(i-1)\tilde{T}}(\omega).$$
    We want to bound $\overline{K}^m(\omega)$ independently of $m$. To do this, denote, for $1< i\lq k$,
    $$A_1=R(\xi)\qtx{and}A_i=A_{i-1}+\sup\big\{R(a)\,|\,a\in\overline{B}(\xi,A_{i-1})\big\}.$$
    By \ref{loccb} of Proposition \ref{pexistloc} we have $\abs{\xi^m_i-\xi^m_{i-1}}\lq R(\xi^m_{i-1})$ for all $1< i\lq k$.
    By induction it follows that $\abs{\xi^m_i-\xi}\lq A_{i-1}$ and $R(\xi^m_i)\lq A_k$ for every $1\lq i\lq k$. Thus,
    writing definition of $K$ from \eqref{defK} and using that $F$ is $\gamma$-Hölder,
    \begin{align*}
        \overline{K}^m(\omega)&=\max_{i=1}^k \big(1+\norm*{F(\xi^m_i)}+H_\gamma(F)R(\xi^m_i)\big)\\
        &\lq 1+\norm{F(\xi)}+2H_\gamma(F)A_k^\gamma\eqcolon \overline{K}(\omega).
    \end{align*}
    This shows \ref{globcc}. Finally, for $p>1/(\gamma\alpha)$, thanks to \eqref{realestvp} and classical rules of the $p$-variation, we have
    \begin{align*}
        V_p(v^m(\omega)) & \lq k^{1-1/p}\sum_{i=1}^k\big(V_p\big(v^m(\omega),{[\smash{(i-1)\tilde{T},i\tilde{T}}[}\big)                                                       \\[-1.5ex] & \hspace{5.75em}+\osc\big(v^m(\omega),[(i-1)\tilde{T},i\tilde{T}]\big)\big)                                                                                          \\[1ex]
                         & \lq k^{2-1/p}\big(C_2(p)\overline{K}(\omega)^\gamma\tilde{T}^{\gamma\beta}+2\overline{K}(\omega)\big),
    \end{align*}
    and then \ref{globce} is shown. This ends the proof of Proposition \ref{pexistglob}.
\end{proof}
\tikzset{
    mid arrow/.style={
            postaction={
                    decorate,
                    decoration={
                            markings,
                            mark=at position 2/3 with {\arrow{>}},
                        }
                }
        }
}
\tikzset{
    mid arrows/.style={
            postaction={
                    decorate,
                    decoration={
                            markings,
                            mark=at position 1/3 with {\arrow{>}},
                            mark=at position 2/3 with {\arrow{>}},
                        }
                }
        }
}
\begin{figure}[h]
    \begin{center}
        \begin{tikzpicture}[x=.75cm,y=.75cm,>=Stealth]
            \coordinate (A) at (1.6913, -2.0356);
            \coordinate (B) at (6.2846, 0.4052);
            \coordinate (C) at (8.3986, -4.1305);
            \coordinate (D) at (6.4575, -8.0511);
            \coordinate (E) at (4.5985, -7.8445);
            \draw[mid arrows,>={Triangle[scale=.6]},line width=1pt] (A) node {$\bullet$} node[left,xshift=-1mm] {$\xi_1$}.. controls (3.421, -0.2291) and (5.0161, -0.575) .. (B) ;
            \draw[mid arrows,>={Triangle[scale=.6]},line width=1pt] (B) node {$\bullet$} node[above,yshift=1mm,xshift=1mm] {$\xi_2$}.. controls (5.9194, -2.8236) and (8.7638, -2.6891) .. (C);
            \draw[mid arrows,>={Triangle[scale=.6]},line width=1pt] (C) node {$\bullet$} node[right,xshift=1mm] {$\xi_3$}.. controls (7.7452, -6.3022) and (6.4575, -8.0511) .. (D) node {$\bullet$} node[below,yshift=-1mm] {$\xi_4$};
            \path[mid arrow,>={Triangle[scale=.6]},draw=black,line width=1pt, dashed]
            (D) .. controls (5.2561, -8.2156) and (5.0279, -8.0186) ..
            (E);
            \draw[>=|,line width=.4pt,>-<,transform canvas={shift={(-0.8136/2.25,1.5311/2.25)}}] (A) -- node[midway, above,sloped] {$\lq R(\xi_1)=A_1$} (B);
            \draw[>=|,line width=.4pt,>-<,transform canvas={shift={(4.5357/6,2.114/6)}}] (B) -- node[midway, above,sloped] {$\lq R(\xi_2)$} (C);
            \draw[>=|,line width=.4pt,>-<,transform canvas={shift={(3.9206/5,-1.9411/5)}}] (C) -- node[midway, below,sloped] {$\lq R(\xi_3)$} (D);
            \draw[shorten <=2mm, shorten >=2mm,line width=.4pt,<->] (A)  -- node[midway, above,sloped] {$\lq A_2$} (C);
            \draw[shorten <=2mm, shorten >=2mm,line width=.4pt,<->] (A)  -- node[midway, above,sloped] {$\lq A_3$} (D);
        \end{tikzpicture}
    \end{center}
    \vspace{-1em}
    \small\caption{Construction of the global approximate solution.}
    \label{figglob}
\end{figure}
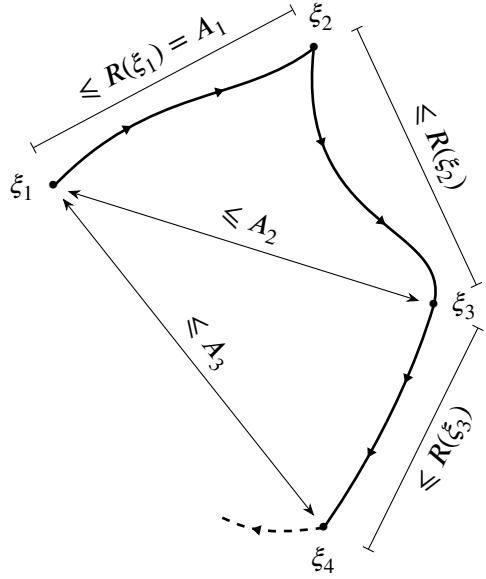
In the rest, $(z^m(\omega),v^m(\omega))_{m,\omega}$ is a
familly of global approximate solutions given by Proposition \ref{pexistglob}.
The next step now is to extract a convergent subsequence of $(v^m(\omega))_m$
for every $\omega\in\Omega$ using a compactness argument.
For us, it is Theorem \ref{thasc}
that makes this extraction possible.
However, because we want, at the end, measurability with respect to $\omega$,
we must perform this extraction in a measurable way.
This is where the results from Section \ref{sectools} come into play.
\begin{npropn}[extraction of limits]\label{pextract}
    For every $\omega=(\xi,x)\in\Omega$, there exist two paths $v(\omega)\in\Vp[1/(\gamma\alpha)+]([0,T]\to\R^{e\times d})$
    and $z(\omega)\in\cl([0,T]\to U)$ with
    \begin{enumtimesrom}
        \item\label{extrca} $z_t(\omega)=\xi+\int_0^t v_s(\omega)\D x_s$ \ for all $t\in[0,T]$.
\item\label{extrcc} $v(\omega)$ is limit of a subsequence of $(v^m(\omega))_m$. Convergence occurs in distance $\dSp[p]$,
        for every $p>1/(\gamma\alpha)$. \item\label{extrcd} $z(\omega)$ is limit of a subsequence of $(z^m(\omega))_m$. Convergence occurs in $\beta$-Hölder norm,
        for every $0<\beta<\alpha$. \item\label{extrce} $z$ and $v$ are measurable from $\Omega$ to their target path spaces.
    \end{enumtimesrom}
\end{npropn}
\begin{proof} In the proof, we will write $\sk$ as a shorthand for $\sk([0,T]\to\R^{e\times d})$.
    The latter set is endowed with distance $\dS$.
    For every $\omega\in\Omega$, consider $\adh_m v^m(\omega)$ the set of subsequential
    limits of $(v^m(\omega))_m$ in distance $\dS$. We look at the map
    $$V:\Omega\ni\omega\longmapsto\adh_m v^m(\omega)\subseteq\sk.$$
    Clearly $V$ is closed-valued. For all $\omega=(\xi,x)\in\Omega$, $m\in\N$ and $t\in[0,T]$,
    $$\abs{v_t^m(\omega)}\lq \abs{v_0^m(\omega)}+\osc(v^m(\omega)).$$
    We have that $v^m_0(\omega)\in F(z^m_0(\omega))=F(\xi)$ thus $\abs{v_0^m(\omega)}\lq\norm{F(\xi)}$.
    Together with $\osc(v^m(\omega))\lq V_p(v^m(\omega))$ for whatever $p>1/(\gamma\alpha)$,
    we obtain that $\norm{v^m(\omega)}_\infty$ is bounded uniformly in $m$.
    Moreover, $(v^m(\omega))_m$ is equioscillating,
    so Theorem \ref{thasc}
    ensures that $(v^m(\omega))_m$ is precompact in uniform norm, and so in distance $\dS$.
In particular the images of $V$ are not empty.
    Let us now show that $V$ is weakly measurable.
    As the space $\sk$ is separable, we can restrict to its open balls.
    So, let $u\in\sk$ and $r>0$.
    Because $(v^m(\omega))_m$ is precompact in distance $\dS$,
    we can write
    $$V^{-1}(B(u,r))=\bigcup_n\limsup_m\left\{\omega:v^m(\omega)\in B\left(u,r-1/n\right)\right\}.$$
    For all $m$, the map $v^m:\Omega\to \sk$ is measurable
    when $\sk$ is endowed with its cylindrical $\sigma$-algebra.
    But we saw in Subsection \ref{ssnotats} that the latter is equal to the Borel $\sigma$-algebra of distance $\dS$.
    Consequently, every set $\left\{\omega:v^m(\omega)\in B\left(u,r-1/n\right)\right\}$
    is measurable in $\Omega$. Hence $V^{-1}(B(u,r))$ is a measurable subset of $\Omega$.
    We can therefore take a measurable selection $v$ of $V$.
    For every $\omega=(\xi,x)\in\Omega$, there is some subsequence $(v^{\phi(m)}(\omega))_m$ that converges to $v(\omega)$ in distance $\dS$. It means that there is a sequence
    $(\lambda_m)_m\subseteq\Lambda$ such that $(v^{\phi(m)}(\omega)\circ\lambda_m)_m$
    converges uniformly to $v(\omega)$. Thus, for any $p>1/(\gamma\alpha)$,
    $$V_p(v(\omega))\lq\liminf_mV_p(v^{\phi(m)}(\omega)\circ\lambda_m)=\liminf_mV_p(v^{\phi(m)}(\omega)).$$
    Because $V_p(v^m(\omega))_m$ is uniformly bounded in $m$, we deduce that $V_p(v(\omega))<\infty$.
    Moreover, applying interpolation inequality \eqref{ineqinterpvp} of Section \ref{sectools} we obtain
    convergence in $\dSp$ of $(v^{\phi(m)}(\omega))_m$ to $v(\omega)$, which shows \ref{extrcc}.
    We set, for all $t\in[0,T]$,
    $$z_t(\omega)=\xi+\int_0^t v_s(\omega)\D x_s.$$
    This ensures \ref{extrca}. By Theorem \ref{contiy}, the sequence $(z^{\phi(m)}(\omega))_m$ converges
    in $\beta$-Hölder norm to $z(\omega)$ for any $0<\beta<\alpha$, which shows \ref{extrcd}.
    For every $t\in[0,T]$, the map $v_t$ is measurable
    because $v$ is measurable. The map $z_t$ is then measurable as a limit
    $$z_t(\omega)=\xi+\lim_n\sum_{[a,b]\in\sigma^n}v_a(\omega)x_{a\wedge t,b\wedge t},$$
    where $(\sigma^n)_n\subseteq\sym_0^T$ is any sequence with $\abs{\sigma^n}\cv[n]0$.
    Hence \ref{extrce} is checked.
\end{proof}
Combining item \ref{globcb} of \ref{pexistglob} with Proposition \ref{pextract}
we can prove Theorem \ref{thydi} for the restricted $\Omega$ defined at the begining
of the subsection.
\begin{proof}[Proof of Lemma \ref{lemrest}]
    Let $\omega\in\Omega$. We show that the limit $(z(\omega),v(\omega))$
    given by \ref{pextract} is solution to \eqref{ydi}.
    Consider the set $\T=\bigcup_m\T_m$ given in Proposition \ref{pexistglob}.
    For $t\in\T$, the bounded sequence $(v^m_t(\omega))_m$
    has a subsequential limit $l$. By \ref{pskcvp}, necessarily $l=v_{t^-}(\omega)$ or $l=v_t(\omega)$.
    But, for $m$ large enough, $t\in\T_m$ so $v^m_t(\omega)\in F(z^m_t(\omega))$.
    Thus, $l\in\Ls_m F(z^m_t(\omega))$. We have \mbox{$z^m_t(\omega)\to z_t(\omega)$}, and,
    since $F$ is continuous, $\Ls_m F(z^m_t(\omega))=F(z_t(\omega))$ so $l\in F(z_t(\omega))$.
    Hence
    \begin{equation}
        l\in \{v_t(\omega),v_{t^-}(\omega)\}\cap F(z_t(\omega))\ne\emptyset.\label{ne}
    \end{equation}
    Moreover if $t\in\{0,T\}$, we have directly $l=v_t(\omega)$ by Skorokhod convergence,
    and thus $v_t(\omega)\in F(z_t(\omega))$.
It remains to consider $t\in{]0,T[}$.
                Take any sequence $(t_n)_n\subseteq{]t,T[}\cap\T$
    converging to $t$. As $v$ is right-continuous,
    $$v_t(\omega)=\lim_n v_{t_n}(\omega)=\lim_n v_{t_n^-}(\omega).$$
    By \eqref{ne} we deduce $v_t(\omega)\in\Li_n F(z_{t_n}(\omega))$.
    Since $z(\omega)$ and $F$ are continuous, we have $\Li_n F(z_{t_n}(\omega))=F(z_t(\omega))$
    and hence $v_t(\omega)\in F(z_t(\omega))$.
\end{proof}
\subsection{Estimates on velocity}\label{sspest}
In the sequel, $S>0$ is a finite time horizon,
and $D^m$ stands for $D^m(S)$.
We consider a path $x:[0,S]\to\R^d$ with $H_\alpha(x)\lq 1$ for some $1/2<\alpha\lq 1$.
This subsection is dedicated to proving the following technical lemma,
which contains all the estimates needed to construct approximate solutions
in Proposition \ref{pexistloc}. This lemma is very close to \cite[Prop.\,2.1]{bailleul_young_2020};
the main difference is that, by allowing the constant $K$ in estimates, and using the additional hypothesis $H_\alpha(x)\lq 1$,
we obtain a less restrictive bound on $S$.
\begin{nlem}\label{techlem1}
    Fix $m\in\N$ and let $v:{[0,S[}\to\R^{e\times d}$ be a path which is constant on
    each interval $[D^m_{k-1},D^m_k[$\, with $k=1,\dots,2^m$.
    Let us denote for all $t\in[0,S]$
    $$z_t=\int_0^t v_s\D x_s=\sum_{k=1}^{2^m}v\left(D^m_{k-1}\right)x\big(D^m_{k-1}\wedge t,D^m_k\wedge t\big).$$
    Suppose that $v$ is bounded by some constant $K\gq 1$ and that,
    for some constant $H$ and all $t\in D^m\ssm\{0,S\}$, we have
    \begin{equation}
        \abs{v_{t}-v_{a(t)}}\lq H\abs{z_{t}-z_{a(t)}}^\gamma.\label{defances}
    \end{equation}
    Pick $1/(1+\gamma)<\beta<\alpha$ and $0<\theta<1$ and suppose also that \begin{equation}\label{bornt}
            S\lq 2^{-\frac{1}{\theta(\alpha-\beta)}} \wedge\left(\frac{H}{1-2^{1-(\gamma+1)\alpha}}\right)^{-\frac{1}{\theta(\alpha-\beta)}}
            \wedge\left(\frac{2H}{1-2^{-\alpha\gamma}}\right)^{-\frac{1}{\gamma(1-\theta)(\alpha-\beta)}}.
    \end{equation}
    Then, for all $n\in\N$ and any two consecutive points $s,t$ of $D^n$,
    \begin{equation}
        \osc\left(v,{[s,t[}\right)\lq K^\gamma 2^{-\gamma\alpha n}S^{\gamma\beta}\ind_{n<m}.\label{estosc}
    \end{equation} Moreover, for any $p>1/(\gamma\alpha)$, we have the estimate
    \begin{equation}
        V_p(v,{[0,S[})\lq \left(\frac{2}{1-2^{1-p\gamma\alpha}}\right)^{1/p}K^\gamma S^{\gamma\beta}.\label{estvp}
    \end{equation}
\end{nlem}
To prove \ref{techlem1} we need the following two intermediate lemmas.
The first one is precisely \cite[Lem.\,2.2]{bailleul_young_2020};
we give its proof here for the sake of completeness. The second lemma
is an enhanced version of \cite[Cor.\,2.3]{bailleul_young_2020}.
\begin{nlem}
    Under the hypotheses of Lemma \ref{techlem1},
    for every $n\lq m$, if $s,t$ are two consecutive points of $D^n$, then
    \begin{equation}\label{intex}
        z_{s,t}=v_sx_{s,t}+\sum_{k=0}^{m-n-1}\sum_{i=0}^{2^k-1}v\big(s^k_i,s^{k+1}_{2i+1}\big)x\big(s^{k+1}_{2i+1},s^k_{i+1}\big)
    \end{equation}
    where $s^k_i=s+iS2^{-n-k}$ for $k\in\N$ and $i=0,\dots,2^k$. See figure \ref{figski}
    for an illustration of the numbers $s^k_i$.
\end{nlem}
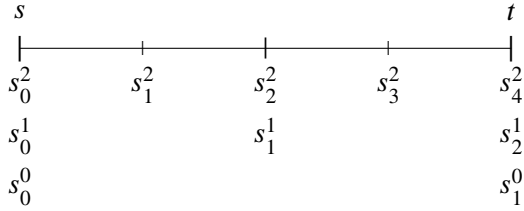
\begin{figure}[h]
    \begin{center}
        \begin{tikzpicture}[x=.5\textwidth]
            \draw (0,0) -- (1,0);
            \draw (0,0) node {$|$} node[yshift=-1.25em] {$s^2_0$} node[yshift=-3em] {$s^1_0$} node[yshift=-4.75em] {$s^0_0$} node[yshift=1.25em] {$s$}
            (1/4,0) node {\tiny$|$} node[yshift=-1.25em] {$s^2_1$}
            (1/2,0) node {\footnotesize$|$} node[yshift=-1.25em] {$s^2_2$} node[yshift=-3em] {$s^1_1$}
            (3/4,0) node {\tiny$|$} node[yshift=-1.25em] {$s^2_3$}
            (1,0) node {$|$} node[yshift=-1.25em] {$s^2_4$} node[yshift=-3em] {$s^1_2$} node[yshift=-4.75em] {$s^0_1$} node[yshift=1.25em] {$t$};
        \end{tikzpicture}
    \end{center}
    \vspace{-1em}
    \small\caption{Values of $s^k_i$ for $k=0,1,2$.}
    \label{figski}
\end{figure}
\begin{proof}
    We prove this lemma by descending induction on $n$. Equation \eqref{intex} is clearly true for $n=m$.
    Suppose that \eqref{intex} is true for some $0<n\lq m$ and let us prove that \eqref{intex} is true for $n-1$.
    Let $s,t$ be two consecutive points of $D^{n-1}$. We can split
    \begin{align*}
        z_{s,t}=z\big(s,s^1_1\big)+z\big(s^1_1,t\big).
    \end{align*}
    Because $s,s^1_1$ and $s^1_1,t$ are consecutive points of $D^n$, we can apply the induction hypothesis,
    and we obtain
    \begin{align*}
        z_{s,t} & =v(s)x\big(s,s^1_1\big)+v\big(s^1_1\big)x\big(s^1_1,t\big)                                                                                             \\
                & \phantom{=\ }+\sum_{k=0}^{m-n-1}\sum_{i=0}^{2^k-1}v\big(s^{k+1}_i,s^{k+2}_{2i+1}\big)x\big(s^{k+2}_{2i+1},s^{k+1}_{i+1}\big)                           \\
& \phantom{=\ }+\sum_{k=0}^{m-n-1}\sum_{i=0}^{2^k-1}v\big(s^{k+1}_{2^k+i},s^{k+2}_{2^{k+1}+2i+1}\big)x\big(s^{k+2}_{2^{k+1}+2i+1},s^{k+1}_{2^k+i+1}\big) \\
                & =v(s)x\big(s,s^1_1\big)+v\big(s^1_1\big)x\big(s^1_1,t\big)
        +\sum_{k=1}^{m-n}\sum_{i=0}^{2^k-1}v\big(s^k_i,s^{k+1}_{2i+1}\big)x\big(s^{k+1}_{2i+1},s^k_{i+1}\big).
    \end{align*}
    We can rewrite the two first terms of the right-hand side as
    \begin{align*}
        v(s)x\big(s,s^1_1\big)+v\big(s^1_1\big)x\big(s^1_1,t\big) & =
        v_sx_{s,t}-v(s)x\big(s^1_1,t\big)+v\big(s^1_1\big)x\big(s^1_1,t\big)                                          \\
                                                                  & =v_sx_{s,t}+v\big(s,s^1_1\big)x\big(s^1_1,t\big).
    \end{align*}
    By replacing this in the previous equality, we obtain \eqref{intex} for $n-1$.
\end{proof}
\begin{nlem}
    Under the hypotheses of Lemma \ref{techlem1},
    for every $n\lq m$, if $s,t$ are two consecutive points of $D^n$, then
    \begin{equation}
        \abs{z_{s,t}}\lq K  2^{-\alpha n} S^{(1-\theta)\alpha+\theta\beta}.\label{estonz}
    \end{equation}
\end{nlem}
\begin{proof}
    For $n\in\N$ we set $B_n=K2^{-\alpha n}S^{(1-\theta)\alpha+\theta\beta}$. We do the proof
    by descending induction on $n$. If $n=m$, then, for $s,t$ two consecutive points of $D^m$,
    \begin{align*}
        \abs{z_{s,t}}=\abs{v_sx_{s,t}} \lq K(S2^{-m})^\alpha =B_mS^{\theta(\alpha-\beta)}\lq B_m.
    \end{align*}
    The last inequality is because $S\lq 1$, due to \eqref{bornt}. Suppose now that \eqref{estonz} is
    true for $n+1,\dots,m$ with $0\lq n<m$ and let us prove that \eqref{estonz} is true
    for $n$. Let $s,t$ two consecutive points of $D^n$. By \eqref{intex} we have
    \begin{align*}
        \abs{z_{s,t}} & \lq\underbrace{\abs{v_sx_{s,t}}}_{\circled{\footnotesize1}}+\sum_{k=0}^{m+n-1}\sum_{i=0}^{2^k-1}\underbrace{\big|v\big(s^k_i,s^{k+1}_{2i+1})\big|\big|x\big(s^{k+1}_{2i+1},s^k_{i+1})\big|}_{\circled{\footnotesize2}}.
    \end{align*}
    For the first term,
    $$\circled{\footnotesize1}\lq K (S2^{-n})^\alpha= B_n S^{\theta(\alpha-\beta)}\lq\tfrac{1}{2}B_n.$$
    The last inequality is due to \eqref{bornt}. For the second term,
    since $s^k_i$ is the ancestor of $s^{k+1}_{2i+1}$, we can use \eqref{defances},
    making appear an increment of $z$, and then apply the induction hypothesis.
    This gives
    \begin{align*}
        \circled{\footnotesize2} & \lq H\big|z\big(s^k_i,s^{k+1}_{2i+1}\big)\big|^\gamma (S2^{-n-k-1})^\alpha \\
                                 & \lq HB_{n+k+1}^\gamma (S2^{-n-k-1})^\alpha                                 \\
                                 & \lq K H2^{-(\gamma+1)\alpha(n+k+1)}S^\alpha.
    \end{align*}
    The last inequality is obtained by writing the definition of $B_{n+k+1}$ and using that $S\lq 1$ and $K\gq 1$.
    Combining those inequalities we get
    \begin{align*}
        \abs{z_{s,t}} & \lq\tfrac{1}{2}B_n+KHS^\alpha2^{-(\gamma+1)\alpha (n+1)}\sum_{k=0}^{m+n-1}2^{k(1-(\gamma+1)\alpha)}    \\
                      & \lq\tfrac{1}{2}B_n+\tfrac{1}{2}B_n S^{\theta(\alpha-\beta)}\frac{H}{1-2^{1-(\gamma+1)\alpha}}\lq B_n.
    \end{align*}
    For the last line we used $(\gamma+1)\alpha>1$ given by the hypothesis on $\alpha$ and $\gamma$.
    This way we have proved that \eqref{estonz} is true for $n$.
\end{proof}
\begin{proof}[Proof of estimate \eqref{estosc} of Lemma \ref{techlem1}]
    The estimate is clear when $n\gq m$ because $v$ is constant on the intervals of $D^m$.
    Fix an integer $0\lq n<m$ and let $s,t$ be two consecutive points of $D^n$.
    Let $u$ be an element of $D^m\cap {[s,t[}$.
    One can see that the set
    $\{0\lq k\lq m-n:a^k(u)=s\}$
    is not empty, and so we denote its minimum by $r$. We then have
    \begin{align*}
        \abs{v_{s,u}} & \lq\sum_{k=1}^r\big|v\big(a^k(u),a^{k-1}(u)\big)\big|           \\
                       & \lq H\sum_{k=1}^r\big|z\big(a^k(u),a^{k-1}(u)\big)\big|^\gamma.
    \end{align*}
    The points $a^k(u)$ and $a^{k-1}(u)$ are consecutive points of $D^l$ for some $n< l\lq m$.
    We can then apply \eqref{estonz} and we obtain
    \begin{align*}
        \abs{v_{s,u}} & \lq H K^\gamma\sum_{l=n+1}^{m}2^{-\gamma\alpha l}S^{\gamma((1-\theta)\alpha+\theta\beta)}                      \\
                       & \lq K^\gamma 2^{-\gamma\alpha n}S^{\gamma\beta}\frac{H}{1-2^{-\gamma\alpha}}S^{\gamma(1-\theta)(\alpha-\beta)} \\
                       & \lq \frac{1}{2}K^\gamma 2^{-\gamma\alpha n}S^{\gamma\beta}.
    \end{align*}
    The last inequality is due to \eqref{bornt}. Therefore, for any $u,w\in D^m\cap{[s,t[}$,
    $$\abs{v_{u,w}}\lq\abs{v_{s,u}}+\abs{v_{s,w}}\lq K^\gamma 2^{-\gamma\alpha n}S^{\gamma\beta}.$$
    As $v$ is constant on the intervals of $D^m$, we get the desired bound on $\osc(v,{[s,t[})$.
\end{proof}
\begin{proof}[Proof of estimate \eqref{estvp} of Lemma \ref{techlem1}]
    Let $(\sigma_k)_{k\lq N}\in\sym({[0,S[})$ a subdivision.
    For every $n\in\N$ consider the set of indices
    $$J_n=\{1\lq k\lq N:\exists\, s,t\text{ consecutive in $D^n$}: s\lq\sigma_{k-1}<\sigma_k<t\}.$$
    For any $k\in J_n$, we can use \eqref{estosc} to get
    $$\abs{v(\sigma_{k-1},\sigma_k)}\lq K^\gamma 2^{-\gamma\alpha n}S^{\gamma\beta}\ind_{n<m}.$$
    Notice that $J_0=\{1,\dots,N\}$ and $J_{n+1}\subseteq J_n$. We then want to bound
    the number of elements of the set $K_n=J_n\ssm J_{n+1}$. To do this,
    one can see that the complementary of $J_n$ is the set
    $$I_n=\{1\lq k\lq N:\exists\, t\in D^n:\sigma_{k-1}<t\lq\sigma_k\}.$$
    The set $I_n$ injects into $D^n\smallsetminus\{0\}$, so its cardinal is less than $2^n$.
    Now, $K_n=I_{n+1}\ssm I_n$ and so $\abs{K_n}\lq\abs{I_{n+1}}\lq 2^{n+1}$. We then have,
    for $p>1/(\gamma\alpha)$,
    \begin{align*}
        V_p[v,\sigma]^p & =\sum_{n=0}^{m-1}\sum_{k\in K_n}\abs{v(\sigma_{k-1},\sigma_k)}^p          \\ & \lq\sum_{n=0}^{m-1}2^{n+1}K^{p\gamma} 2^{-p\gamma\alpha n}S^{p\gamma\beta} \\
                        & \lq\frac{1}{1-2^{1-p\gamma\alpha}}2K^{p\gamma}S^{p\gamma\beta}.            
    \end{align*}
    This concludes the proof of \eqref{estvp}.
\end{proof}
\subsection{General case}\label{sspglobgen}
In this subsection we adopt the assumptions of Theorem \ref{thydi}.
For $S>0$ we define
\begin{align*}
    \tstc([0,S])&=\{x:[0,S]\to\R^d:H_{1/q}(x)\lq 1\},\\
    \F([0,S])&=\{x:[0,T]\to[0,S]\text{ continuous nondecreasing surjective}\}.
\end{align*}
We will derive the proof of Theorem \ref{thydi} from Lemma \ref{lemrest} by splitting a path $x\in\vpc[q]([0,T]\to\R^d)$
into a composition of a $1/q$-Hölder path and a nondecreasing path, like in \cite[\S\,3]{chistyakov_maps_1998}. More precisely, we set, for every integer $n\gq1$,
$$\G_n=\big\{x\in\vpc[q]([0,T]\to\R^d):(n-1)T\lq V_q(x)^q<nT\big\}.$$
Therefore the set $\vpc[q]([0,T]\to\R^d)$ is just the disjoint union of $(\G_n)_n$.
For every $x\in\G_n$ we set
\begin{align*}
    \forall t\in[0,T],\ \ \Phi x(t)&=\frac{nT}{V_q(x)^q}V_q(x,[0,t])^q,\\
    \forall c\in[0,nT],\ \ \Psi x(c)&=x((\Phi x)^{-1}(c)).
\end{align*}
The map $\Psi x$ is well defined because, for $c\in[0,nT]$, if $s\lq t\in(\Phi x)^{-1}(c)$,
$$\abs{x_t-x_s}^q\lq V_q(x,[s,t])^q\lq \Phi x(t)-\Phi x(s)=0.$$
The map $[0,T]\ni t\mapsto V_q(x,[0,t])$ is continuous for $x\in\vpc[q]([0,T]\to\R^d)$,
see \cite[Th.\,4.2]{chistyakov_maps_1998}.
Thus, one checks that $\Phi x\in\F([0,nT])$, $\Psi x\in\tstc([0,nT])$ and $x=\Psi x\circ\Phi x$.
Applying Lemma \ref{lemrest} to $\Psi x$ we get a measurable map $$(z^n,v^n):U\times\tstc([0,nT])\to\cl[1/q]([0,nT]\to U)\times \Vp[q/\gamma+]([0,nT]\to\R^{e\times d}).$$
We then set $\Omega_n=U\times\G_n$ and, for $\omega=(\xi,x)\in \Omega_n$,
$$(z,v)(\omega)=(z^n,v^n)(\xi,\Psi x)\circ\Phi x.$$
This way, $(z,v)$ is a map from $\Omega$ to $(\vpc[q]([0,T]\to U))\times \Vp[q/\gamma+]([0,T]\to\R^{e\times d})$.
By \ref{pyicv}, the path $z(\omega)$ is indeed a solution to \eqref{ydi} with velocity $v(\omega)$.
It remains to check that $(z,v)$ is measurable. The two following lemmas will be useful.
\begin{nlem}\label{lemvpmes}
    The map $x\mapsto V_q(x)$ is measurable on $\vpc[q]([0,T]\to\R^d)$.
\end{nlem}
\begin{proof}
    For all $x\in \vpc[q]([0,T]\to\R^d)$ we have by continuity of $x$
    $$V_q(x)=\sup_{\substack{(\sigma_k)_{k\lq n}\in\sym([0,T])\\\sigma\subseteq\Q}}\left(\sum_{k=1}^n\abs{x(\sigma_{k-1},\sigma_k)}^q\right)^{1/q},$$
    hence the application $x\mapsto V_q(x)$ is a countable supremum of measurable functions so is measurable itself.
\end{proof}
\begin{nlem}\label{lemcompmes}
    For $S>0$ and $Y$ a separable metric space, the map
    $$\Theta:\sk([0,S]\to Y)\times[0,S]^{[0,T]}\ni(x,\lambda)\longmapsto x\circ\lambda\in Y^{[0,T]}$$
    is measurable --- with respect to cylindrical $\sigma$-algebras.
\end{nlem}
\begin{proof}
    As $Y$ is separable we can restrict to its open balls. Let $y\in Y$ and $r>0$. By right-continuity
    of $x$ we can write, for all $t\in[0,T]$,
    \begin{align*}
        \Theta_t^{-1}(B(y,r))&=\{x:x(S)\in B(y,r)\}\times \{\lambda:\lambda(t)=S\}\,\cup\\
        &\hspace{-4em}\bigcup_n\bigcap_m\bigcup_{a\in[0,S]\cap\Q}\{x:x(a)\in B(y,r-1/n)\}\times \{\lambda:a-1/m\lq\lambda(t)\lq a\}.
    \end{align*}
    Hence $\Theta_t^{-1}(B(y,r))$ is in $\B(Y)^{\otimes[0,S]}\otimes\B([0,S])^{\otimes[0,T]}$, so $\Theta$ is measurable.
\end{proof}
A consequence of Lemma \ref{lemvpmes} is that every set $\G_n$ is measurable. Now we can show that the maps $\Phi$ and $\Psi$ are measurable.
It suffices to show measurability on every $\G_n$.
\begin{nlem}
    The applications $\Phi:\G_n\to\F([0,nT])$ and $\Psi:\G_n\to\tstc([0,nT])$ are measurable.
\end{nlem}
\begin{proof}
    By Lemma \ref{lemvpmes} it is clear that $\Phi$ is measurable. For every $x\in\G_n$, we have $\Psi x=\Theta(x,\Gamma x)$,
    where $$\Gamma x:[0,nT]\ni c\longmapsto \max(\Phi x)^{-1}(c)\in[0,T]$$
    and $\Theta$ is defined in Lemma \ref{lemcompmes}.
    As $\Theta$ is measurable we just have to check that $\Gamma$ is measurable. For $c\in[0,nT]$ and $r\in[0,T]$,
    \begin{align*}
        \Gamma x(c)\lq r\,&\Longleftrightarrow\, c\notin \Phi x({]r,T]})\\
        &\Longleftrightarrow\, \Phi x({]r,T]})<c\text{ \,or \,}\Phi x({]r,T]})>c\\
        &\Longleftrightarrow\, \Phi x(T)<c\text{ \,or \,}\forall t\in{]r,T]}\cap\Q, \Phi x(t)>c.
    \end{align*}
    The second line is due to intermediate value theorem, as $\Phi x$ is continuous. The last line is
    because $\Phi x$ is nondecreasing. Since $\Phi$ is measurable we see that $\Gamma$ also is.
\end{proof}
Therefore, the map $(z,v)$ appears as a composition of multiple measurable maps, so is measurable itself,
which concludes the proof of \ref{thydi} in the general case.
\section{Other proofs}\label{secproof2}
\subsection{Local existence result}

The key of the proofs in this subsection is the following extension lemma.
\begin{nlem}[extension of Hölder set-valued maps]\label{lemext}
    Let $A$ be a closed convex subset of an Hilbert space $X$.
    Let $F$ be a map defined on $A$ and valued in the nonempty subsets of a metric space $Y$.
    If $F$ is $\gamma$-Hölder for some $0<\gamma\lq 1$, then $F$ extends to
    a $\gamma$-Hölder set-valued map defined on $X$ with the same Hölder constant and the same images as $F$.
\end{nlem}
\begin{proof}
    For every $x\in X$ define $\tilde{F}(x)=F(p(x))$
    where $p$ denotes the orthogonal projection on $A$.
    Hence $\tilde{F}$ coincides with $F$ on $A$ and,
    because $p$ is $1$-Lipschitz, $\tilde{F}$ is $\gamma$-Hölder
    with the same Hölder constant as $F$.
\end{proof}
\begin{proof}[Proof of Theorem \ref{thydiloc}]
    Here $F$ is valued in $\kl(\R^{e\times d})$ and locally $\gamma$-Hölder for some $q-1<\gamma\lq 1$.
    Write $U$ as a countable union of open balls $B_n$ and define, for every $n$,
    $$C_n=B_n\smallsetminus\bigcup_{m<n} B_m\qtx{and}\Omega_n=C_n\times\vpc[q]([0,T]\to\R^d),$$
    so that $\Omega_n$ is measurable and $\Omega$ is the disjoint union of $(\Omega_n)_n$.
    By assumption $F$ is $\gamma$-Hölder on $\overline{B_n}$, so let $F_n$ be a $\gamma$-Hölder map
    valued in $\kl(\R^{e\times d})$ and defined on $\R^e$ that coincides with $F$ on $\overline{B_n}$,
    see Lemma \ref{lemext}.
    Theorem \ref{thydi} applied with vector field $F_n$ gives us a pair $(z^n,v^n)$.
    For $\omega=(\xi,x)\in\Omega$, there is only one $n$ such that $\omega\in \Omega_n$.
    We then set
    $$T_\ast(\omega)=\inf\,\{t\in[0,T]:z^n_t(\omega)\notin\overline{B_n}\}$$
    and, for every $t\in[0,T_\ast(\omega)]$,
    $$(z_t(\omega),v_t(\omega))=(z^n_t(\omega),v^n_t(\omega))$$
    Since $z^n(\omega)$ is continuous, we have $T_\ast(\omega)>0$,
    and $z^n(\omega)$ stays in $\overline{B_n}$ on $[0,T_\ast(\omega)]$.
    As $F_n=F$ on $\overline{B_n}$, the pair $(z(\omega),v(\omega))$ is hence a solution to \eqref{ydi} started at $\xi$ and whose velocity is
    in $\Vp[q/\gamma+]([0,T_\ast(\omega)]\to\R^{e\times d})$.
    For every $0<A\lq T$,
    $$T_\ast(\omega)\gq A\Longleftrightarrow\omega\in\bigcup_n\left(\Omega_n\cap \bigcap_{t\in[0,A]\cap\Q}(z^n_t)^{-1}\big(\overline{B_n}\big)\right).$$
    Hence $T_\ast$ is measurable, then one easily checks that the map $(z,v)\ind_{[0,T_\ast]}$ is measurable, which ends the proof of the theorem.
\end{proof}
The construction of the next proof is a bit more complicated. It requires the following measurability lemma.
\begin{nlem}\label{lemmes}
    Consider a map $y:\Omega\to\sk([0,T]\to Y)$ with $\Omega$ a measurable space
    and $Y$ a separable metric space. Then, $y$ is measurable
    if and only if $y$ is measurable from the product space $[0,T]\times \Omega$ to $Y$.
\end{nlem}
\begin{proof}
    Suppose $y$ is measurable. For $u\in Y$ and $r>0$, we can write
    $$\{(\omega,t):y_t(\omega)\in B(u,r)\}=\bigcup_n\bigcap_k\bigcup_{t\in E}y_t^{-1}(B(u,r-1/n))\times B(t,1/k),$$
    where $E$ is a dense subset of $[0,T]$ containing $T$. As $Y$ is separable this suffices to prove
    that $y$ is measurable from $[0,T]\times \Omega$ to $Y$. The converse is obvious.
\end{proof}
\begin{proof}[Proof of Theorem \ref{thydimax}]
    Here $U$ is convex and $F$ is valued in $\kl(\R^{e\times d})$ and locally $\gamma$-Hölder
    for some $q-1<\gamma\lq 1$.
    For every $n\in\N$, we set $$K_n=\left\{a\in U:\dist(a,U\cp)\gq1/n\right\}\cap\overline{B}(0,n)$$
    so that $U$ is the nondecreasing union of $(K_n)_n$.
    By assumption, $K_n$ is compact convex and $F$ is $\gamma$-Hölder on $K_n$.
    Let $F_n$ be a $\gamma$-Hölder
    map valued in $\kl(\R^{e\times d})$ and defined on $\R^e$ that coincides with $F$ on $K_n$,
    see Lemma \ref{lemext}.
    Theorem \ref{thydi} applied with vector field $F_n$ gives us a pair $(y^n,u^n)$.
    For $\omega=(\xi,x)\in\Omega$ define $T_0(\omega)=0$, $$z^0(\omega):[0,T]\ni t\mapsto\xi\qtx{and}v^0(\omega):[0,T]\ni t\mapsto0.$$
    Suppose now that for some $n\in\N$ we have constructed $T_n(\omega)\in[0,T]$ and two paths
    $z^n(\omega)\in\vpc[q]([0,T]\to K_n)$ and $v^n(\omega)\in\Vp[q/\gamma+]([0,T]\to\R^{e\times d})$, such that
    \begin{enumtimes}
        \item\label{maxpx} $(z^n(\omega),v^n(\omega))$ is a solution to \eqref{ydi} started at $\xi$ on $[0,T_n(\omega)]$.
        \item\label{maxpa} $T_n(\omega)\gq T_{n-1}(\omega)$ with strict inequality if $T_{n-1}(\omega)<T$.
        \item\label{maxpb} $z^n(\omega)=z^{n-1}(\omega)$ on $[0,T_{n-1}(\omega)]$ and $v^n(\omega)=v^{n-1}(\omega)$ on ${[0,T_{n-1}(\omega)[}$.
        \item\label{maxpc} $z^n_{T_n(\omega)}(\omega)\in \partial K_n$ if $T_n(\omega)<T$.
        \item\label{maxpd} The maps $T_n$, $z^n$ and $v^n$ are measurable.
    \end{enumtimes}
    If $T_n(\omega)=T$, define $T_{n+1}(\omega)=T$ and $(z^{n+1},v^{n+1})(\omega)=(z^n,v^n)(\omega)$.
    Else, set $z^{n+1}(\omega)=z^n(\omega)$ on $[0,T_n(\omega)]$ and $v^{n+1}(\omega)=v^n(\omega)$ on ${[0,T_n(\omega)[}$,
    which ensures \ref{maxpb} at rank $n+1$, and, for every $t\in[T_n(\omega),T]$,
    $$\big(z^{n+1}_t(\omega),v^{n+1}_t(\omega)\big)=\big(y^{n+1}_{t-T_n(\omega)},u^{n+1}_{t-T_n(\omega)}\big)\big(z^n_{T_n(\omega)}(\omega),x\big)$$
    so that $z^{n+1}(\omega)$ and $v^{n+1}(\omega)$ are defined on $[0,T]$.
    By Lemma \ref{lemmes} the maps $z^{n+1}$ and $v^{n+1}$ are measurable.
    Define then
    $$T_{n+1}(\omega)=\inf\big\{t\in[T_n(\omega),T]:z^{n+1}_t(\omega)\notin K_{n+1}\big\}.$$
    Because $z^{n+1}(\omega)$ is continuous, it is valued in $K_{n+1}$ on $[0,T_{n+1}(\omega)]$ and we have \ref{maxpa},
    \ref{maxpx} and \ref{maxpc} at rank $n+1$. For every $0<A\lq T$ one can see that
    $$\{\omega:T_{n+1}(\omega)\gq A\}=\{\omega:T_n(\omega)\gq A\}\cup\bigcap_{t\in[0,A]\cap\Q}\big\{\omega:z^{n+1}_{T_n(\omega)\vee t}(\omega)\in K_{n+1}\big\}.$$
    Using again Lemma \ref{lemmes} we get that $T_{n+1}$ is measurable and so \ref{maxpd} at rank $n+1$.
    Having run this induction we set $T_\ast(\omega)=\sup_n T_n(\omega)$.
    The map $T_\ast$ is measurable as a supremum.
    If $T_n(\omega)=T$ for some $n$,
    set $I(\omega)=[0,T]$ and $(z,v)(\omega)=(z^n,v^n)(\omega)$,
    otherwise set $I(\omega)={[0,T_\ast(\omega)[}$ and define
    $(z_t(\omega),v_t(\omega))=(z^n_t(\omega),v^n_t(\omega))$,
    where $n$ is such that $t< T_n(\omega)$, for all $t\in I(\omega)$. This way, $(z(\omega),v(\omega))$ is a solution
    to \eqref{ydi} started at $\xi$ whose velocity is in $\Vploc[q/\gamma+](I(\omega)\to\R^{e\times d})$,
    and one checks that $(z,v)\ind_I$ is measurable. In the case where $I(\omega)\ne[0,T]$, observe that, for every $n$, we have $z_{T_n(\omega)}(\omega)\in\partial K_n$ and
    $$\partial K_n \subseteq \left\{a\in U:\dist(a,U\cp)=1/n\right\}\cup S(0,n),$$
    hence, every convergent subsequence of $(z_{T_n(\omega)}(\omega))_n$ has a limit in the closed set $U\cp$,
    which means that $z_{I(\omega)}(\omega)$ is not precompact in $U$.
    This concludes the proof.
\end{proof}
\subsection{Selection result}
The proof of Theorem \ref{thsel} follows the same steps as Theorem \ref{thydi}.
We therefore combine the arguments into a single, more concise proof.
Notations and theorems from the introduction of Section \ref{secproof1} are used.
In this subsection $D^m$ will stand for $D^m(T)$.
\begin{proof}[Proof of Theorem \ref{thsel}]
    Let $m\in\N$. We set, for all $\omega\in\Omega$ and $t\in[0,D^m_1[$,
    $$f^m_t(\omega)=\omega.$$
    Suppose now that for all $\omega\in\Omega$ we have constructed the path $f^m(\omega)$
    on the interval ${[0,\tau[}$, with $\tau\coloneqq D^m_k$
    for some $1\lq k< 2^m$, and that
    \begin{enumtimes}
        \item\label{selcpb} $f^m_t(\omega)\in F(t)$ \ for every $t\in D^m\cap{[0,\tau[}$.
        \item\label{selcpc} $f^m(\omega)$ is constant on every interval $[D^m_{i-1},D^m_i[$ with $i=1,\dots,k$.
        \item\label{selcpe} For all $t\in D^m\cap{]0,\tau[}$, we have
        $$\abs{f^m_{t}(\omega)-f^m_{a(t)}(\omega)}\lq H_\gamma(F)\big(T2^{-M(t)}\big)^\gamma.$$
        \item\label{selcpf} Every map
        $\Omega\ni\omega\mapsto f^m_t(\omega)$,
        with $t\in[0,D^m_k[$, is measurable.
    \end{enumtimes}
    As $F$ is $\gamma$-Hölder we have
    $$h(F(\tau),F(a(\tau)))\lq H_\gamma(F)\abs{\tau-a(\tau)}^\gamma=H_\gamma(F)\big(T2^{-M(\tau)}\big)^\gamma\eqcolon r.$$
    On the other hand, we have $f^m_{a(\tau)}(\omega)\in F(a(\tau))$ by \ref{selcpb}. Then, because $F$ is closed-valued,
    the set
    $$G(\omega)\coloneqq\overline{B}\big(f^m_{a(\tau)}(\omega),r\big)\cap F(\tau)$$
    is nonempty and closed. To show that $G$ is measurable, we write, for a closed subset $A$ of $\R^{e\times d}$,
    $$G^{-1}(A)=\big(f^m_{a(\tau)}\big)^{-1}(F(\tau)\cap A+\overline{B}(0,r)),$$
    and conclude with \ref{selcpf}.
    Therefore we can take a measurable selection of $G$ which we denote by $\phi$.
    We then set, for all $\omega\in\Omega$ and $t\in[D^m_k,D^m_{k+1}[$,
    $$f^m_t(\omega)=\phi(\omega).$$
    Propositions \ref{selcpb}, \ref{selcpc}, \ref{selcpe} and \ref{selcpf} are clearly true at rank $k+1$.
    By running this induction, we build the path $f^m(\omega)$ on ${[0,T[}$.
    To finish we set $f^m_T(\omega)=f^m_{T^-}(\omega)$.
    By the same proof as estimate \eqref{estosc} of Lemma \ref{techlem1},
    we have, for any $\omega\in\Omega$ and two consecutive points $s,t$ of $D^n$, with $n\in\N$,
    \begin{equation}\label{estselosc}
        \osc(f^m(\omega),{[s,t[})\lq \frac{2H_\gamma(F)T^\gamma}{1-2^{-\gamma}}2^{-n\gamma}\ind_{n<m}.
    \end{equation}
    Adapting the proof of estimate \eqref{estvp} of Lemma \ref{techlem1} we obtain also
    \begin{equation}\label{estselvp}
        V_p(f^m(\omega))\lq \left(\frac{2}{1-2^{1-p\gamma}}\right)^{1/p}\frac{2H_\gamma(F)T^\gamma}{1-2^{-\gamma}}
    \end{equation}
    for all $p>1/\gamma$. In the following we write $\sk$ as a shorthand for $\sk([0,T]\to\R^d)$.
    The latter set is endowed with distance $\dS$.
    Consider the set $\adh_m f^m(\omega)$ of subsequential values of $(f^m(\omega))_m$ in distance $\dS$. We look at the map
    $$V:\Omega\ni\omega\longmapsto\adh_m f^m(\omega)\subseteq\sk.$$
    Clearly $V$ is closed valued. By \eqref{estselosc} the familly of paths $(f^m(\omega))_m$ is equioscillating.
    Moreover, $\norm{f^m(\omega)}_\infty\lq \norm{F}_\infty$ for every $m\in\N$,
    hence by Theorem \ref{thasc} the sequence $(f^m(\omega))_m$ is precompact in uniform norm,
    and so in distance $\dS$. In particular the images of $V$ are not empty.
    As in the proof of Theorem \ref{thydi},
    we show that $V$ is weakly measurable by writing
    $$V^{-1}(B(u,r))=\bigcup_n\limsup_m\left\{\omega:f^m(\omega)\in B\left(u,r-1/n\right)\right\}.$$
    where $B(u,r)$ the open ball of center $u$ and radius $r$ in $\sk$.
    We can therefore take a measurable selection $f$ of $V$.
    We prove that $f(\omega)$ is a selection of $F$ using right-continuity of $f(\omega)$ on $[0,T]$,
    in the same way as in proof of Theorem \ref{thydi}.
    By \eqref{estselvp} the path $f(\omega)$ is of bounded $p$-variation for every $p>1/\gamma$ and we have the estimate \eqref{estselvpor}.
\end{proof}
\pagebreak \printrefs \apx
\numberwithin{thm}{chapter}
\chapter{Basics on Young integral}\label{apxa}
To put it briefly, Young integral is Riemann-Stieltjes integral
when the driving signal is not of bounded variation. In this appendix
we summarize how Young integral works, following the pedagogical exposition by
Friz and Zhang \cite{friz_differential_2018}.
Let $T>0$ be a finite time horizon and $d,e\gq 1$ be integers.
\begin{nthm}[Young integral]\label{pyi}
    Considering two paths $x\in\vpc([0,T]\to\R^d)$ and $y\in\vp[q]([0,T]\to\R^{e\times d})$,
    with $p,q\gq 1$, we denote, for a subdivision $(\sigma_k)_{k\lq n}\in\sym_0^T$,
    $$y\D x(\sigma)=\sum_{k=1}^n y(\sigma_{k-1})(x(\sigma_k)-x(\sigma_{k-1})).$$
    If $1/p+1/q>1$, then $y\D x(\sigma)$ converges as $\abs{\sigma}\to0$.
    Its limit is called the Young integral of $y$ against $x$ on $[0,T]$ and is denoted $$\int_0^T y\D x=\int_0^T y_s\D x_s.$$
    For any $s\lq t\in[0,T]$ we have the Young-Love estimate
    $$\abs*{\int_s^t y\D x-y_sx_{s,t}}\lq \frac{1}{1-2^{1-1/p-1/q}}V_q(y,{[s,t[})V_p(x,[s,t]).$$
    Furthermore, the Young integral of $y$ against $x$ is the only $\R^d$-valued function $\mu$
    defined on the triangle $\Delta=\{(s,t)\,|\,s\lq t\in[0,T]\}$ satisfying
    \begin{gather*}
        \mu(s,t)=\mu(s,u)+\mu(u,t),\\
        \abs{\mu(s,t)-y_sx_{s,t}}\lq C V_q(y,[s,t])V_p(x,[s,t]),
    \end{gather*}
    for all $s\lq u\lq t$ in a dense subset of $[0,T]$ and some constant $C$.
\end{nthm}
The assumption that $x$ is continuous in Theorem \ref{pyi}
leads to several simplifications compared with \cite[Prop.\,2.4]{friz_differential_2018}, from which this theorem is derived.
In particular, when the driving signal $x$ is continuous,
one may modify the values of $y$ on any co-dense subset of $[0,T]$ without
affecting the value of the integral of $y$ against $x$, which is useful
when gluing solutions of a Young differential inclusion.
\bskip

A consequence of the Young-Love estimate is the two following propositions:
\begin{npropn}\label{pyivp}
    Given $x\in\vpc([0,T]\to\R^d)$ and $y\in\vp[q]([0,T]\to\R^{e\times d})$ with $1/p+1/q>1$,
    the indefinite Young integral of $y$ against $x$ is of bounded $p$-variation with
    \begin{equation}\label{eq10}
            \normVp{\int_0^\cdot y\D x}\lq \CY(p,q)\normVp[q]{y}V_p(x),
    \end{equation}
    for some constant $\CY(p,q)$ depending on $p$ and $q$.
    As $x$ is continuous, the mentioned integral is also continuous by \eqref{eq10}.
\end{npropn}
\begin{npropn}\label{pyiha}
    Given $x\in\cl([0,T]\to\R^d)$ and $y\in\vp[q]([0,T]\to\R^{e\times d})$ with $\alpha+1/q>1$,
    the indefinite Young integral of $y$ against $x$ is $\alpha$-Hölder with
    \begin{equation}\label{eq11}
            H_\alpha\left(\int_0^\cdot y\D x\right)\lq \CY(1/\alpha,q)\normVp[q]{y}H_\alpha(x).
    \end{equation}
\end{npropn}
\begin{nrk}
    In the present paper, the constant $\CY(p,q)$ is defined so that estimates \eqref{eq10} and \eqref{eq11}
    are true at the same time. Therefore, it is not necessarily minimal.
\end{nrk}
The Young integral is compatible with the two basic operations that are bilinearity and change of time:
\begin{npropn}[bilinearity]\label{pyibil}
If $p,q\gq1$ are such that $1/p+1/q>1$, the map
    $$\fun*{\vpc([0,T]\to\R^d)\times \vp[q]([0,T]\to\R^{e\times d})}{\vpc([0,T]\to\R^e)}{(x,y)}{\int_0^Ty\D x}$$
    is bilinear and continuous under usual norms.
\end{npropn}
\begin{npropn}[change of time]\label{pyicv}
    Consider a continuous nondecreasing
    change of time $\lambda:[0,S]\to[0,T]$. Given any two paths $x\in\vpc([0,T]\to\R^d)$ and $y\in\vp[q]([0,T]\to\R^{e\times d})$ with $1/p+1/q>1$,
    we have for any $s\lq t\in[0,S]$,
    $$\int_s^t(y\circ\lambda)\D(x\circ\lambda)=\int_{\lambda(s)}^{\lambda(t)}y\D x.$$
\end{npropn}
\printrefs
\end{document}